\documentclass[12pt]{article}

\usepackage{amssymb,amsmath,amsfonts,amssymb}
-\marginparwidth \oddsidemargin -4mm \evensidemargin
\oddsidemargin
\usepackage{amssymb}
\advance\hoffset by 5mm

\def\@abssec#1{\vspace{.05in}\footnotesize \parindent .2in
{\bf #1. }\ignorespaces}

\newtheorem{theorem}{Theorem}[section]

\newtheorem{lemma}[theorem]{Lemma}
\newtheorem{proposition}[theorem]{Proposition}
\newtheorem{corollary}[theorem]{Corollary}

\def \Rm {\mathbb R}

\def \Sm {\mathbb S}

\title{ Blow up and regularity for fractal Burgers equation}
\author{Alexander Kiselev\thanks{Department of
Mathematics, University of Wisconsin, Madison, WI 53706, USA;
e-mail: kiselev@math.wisc.edu } \and Fedor
Nazarov\thanks{Department of Mathematics, University of Wisconsin,
Madison, WI 53706, USA; e-mail: nazarov@math.wisc.edu} \and Roman
Shterenberg\thanks{Department of Mathematics, University of
Alabama, Birmingham, AL 35294, USA; e-mail: shterenb@math.uab.edu
}}

\begin{document}

\maketitle

\begin{abstract}
The paper is a comprehensive study of the existence, uniqueness,
blow up and regularity properties of solutions of the Burgers
equation with fractional dissipation. We prove existence of the
finite time blow up for the power of Laplacian $\alpha < 1/2,$ and
global existence as well as analyticity of solution for $\alpha
\geq 1/2.$ We also prove the existence of solutions with very
rough initial data $u_0 \in L^p,$ $1 < p < \infty.$ Many of the
results can be extended to a more general class of equations,
including the surface quasi-geostrophic equation.
\end{abstract}

\section{Introduction}

The purpose of this paper is to present several results on
Burgers equation with fractional dissipation
\begin{equation}\label{bur1}
u_t= u  u_x -(-\Delta)^{\alpha}u, \,\,\,u(x,0)=u_0(x).
\end{equation}
We will consider \eqref{bur1} on the circle $\Sm^1.$ Equivalently,
one can consider \eqref{bur1} on the real line with periodic
initial data $u_0(x)$.

The Burgers equation with $\alpha =0$ and $\alpha =1$ has received
an extensive amount of attention since the studies by Burgers in
the 1940s (and it has been considered even earlier by Beteman
\cite{Beteman} and Forsyth \cite{Forsyth}, pp 97--102). If $\alpha
=0,$ the equation is perhaps the most basic example of a PDE
evolution leading to shocks; if $\alpha =1,$ it provides an
accessible model for studying the interaction between nonlinear
and dissipative phenomena.

The Burgers equation can also be viewed as the simplest in the
family of partial differential equations modeling the Euler and
Navier-Stokes equation nonlinearity. 
Recently, there has been increased interest in models involving
fractional dissipation, in particular Navier-Stokes (see
\cite{KP1}) and surface quasi-geostrophic equations (see e.g.
\cite{CCW,CafVas,KNV,Wu} for further references). Fractional
dissipation also appears naturally in certain combustion models
\cite{Matalon}. 
Our goal is to present here, in the most accessible framework of
the Burgers equation, results and techniques that in some cases
apply (with relatively straightforward adjustments) to a wider
class of equations including in particular surface
quasi-geostrophic. Among the results we prove are the global
existence of solutions for $\alpha = 1/2$ (more generally $\alpha
\geq 1/2$), space analyticity of solutions for $\alpha \geq 1/2,$
the existence of solutions with very rough initial data, as well
as blow up in finite time for $\alpha <1/2$. 

Let us now describe in more detail some of the results that we
prove. We denote by $W^s_p$, $s\in{\mathbb R}$, $1\leq
p\leq\infty$ the standard Sobolev spaces. If $p=2$ we use notation
$H^s$, and denote by $\|\cdot\|_s$ the norm in $H^s.$ Without loss
of generality, we will consider equation \eqref{bur1} on the
subspace of mean zero functions. This subspace is preserved by
evolution and contains all non-trivial dynamics. The advantage of
this subspace is that the $H^s$ norm dominates the $L^2$ norm, and
this simplifies the estimates.
 We will also always assume
that the initial data (and so the solution) are real valued.

The case $\alpha >1/2$ is subcritical, and smooth solutions exist
globally. This is a fairly simple fact to prove, using the maximum
principle control of the $L^\infty$ norm for $\alpha \leq 1,$ and
straightforward estimates for $\alpha >1$ (see \cite{Resnick,CC}
for the quasi-geostrophic case and $1/2 \leq \alpha \leq 1$; the
argument transfers to the dissipative Burgers equation without
significant changes). One can also use the scheme of the proof
suggested here for the critical case $\alpha=1/2$. That is why we
just state the result for $\alpha>1/2$ without proof.

\begin{theorem}\label{subcritical}
Assume that $\alpha>1/2,$ and the initial data $u_0(x)$ belongs to
$H^s,$ $s>3/2-2\alpha,$ $s \geq 0.$ Then there exists a unique
global solution of the equation \eqref{bur1} $u(x,t)$ that belongs
to $C([0,\infty),H^s).$ Moreover, this solution is real analytic
in $x$ for $t>0$.
\end{theorem}

We have been unable to find the analyticity claim in the existent
literature, but again the proof is parallel to that for the
$\alpha=1/2$ case, which we will do in detail. In what follows, we
will consider mostly $0<\alpha\leq1/2$. Our first result concerns
the local existence and uniqueness of classical solutions to
Burgers equation with initial data in Sobolev spaces.

Let us denote $C_w([0,T],L^2)$ the class of solutions that are
weakly continuous as functions with values in $L^2.$ \\

\noindent \it Definition. \rm We say that $u\in L^2([0,T],L^2)
\cap C_w([0,T],L^2)$ and such that $du/dt \in L^1([0,T],H^{-1})$
is a weak solution of \eqref{bur1} for $t\in(0,T)$ if for any
smooth periodic function $\varphi(x)$ we have
\begin{equation}\label{burweak1}
(u,\varphi)_t=-\frac12(u^2,\varphi_x)
-(u,(-\Delta)^{\alpha}\varphi),\ \ \hbox{a.e.}\ t\in(0,T),\ \
(u,\varphi)(0)=(u_0,\varphi).
\end{equation}
Note that then $(u,\varphi)(t)$ is absolutely continuous and
\begin{equation}\label{burweak}
(u,\varphi)(t)-(u_0,\varphi)=\int\limits_0^t\left(-\frac12(u^2,\varphi_x)
-(u,(-\Delta)^{\alpha}\varphi)\right)ds
\end{equation}
for any $t \in [0,T].$

\begin{theorem}\label{thm1}
Assume that $0<\alpha\leq1/2,$ and the initial data $u_0(x) \in
H^s,$ $s>3/2-2\alpha$. Then there exists $T=T(\alpha,\|u_0\|_s)>0$
such that there exists a weak solution of the equation
\eqref{bur1} $u(x,t)$ satisfying $u(x,t) \in C([0,T],H^s) \cap
L^2([0,T],H^{s+\alpha}).$ Moreover, $u$ can be chosen to satisfy
$u(x,t) \in C^{\infty}$ for $0<t<T.$ If $v$ is another weak
solution of \eqref{bur1} with initial data $u_0$ such that $v\in
C([0,T],L^2)\cap L^{3/2\delta}([0,T],H^{\delta})$ with some
$\delta\in(1/2,1]$, then $v$ coincides with $u$.
\end{theorem}

\noindent \it Remark \rm  1. In particular, it follows that the
solution $u$ in Theorem~\ref{thm1} solves \eqref{bur1} in
classical sense for
every $t>0.$ \\
\noindent \it Remark 2. \rm It is clear from Theorem~\ref{thm1}
that the solution is unique
in the $C([0,T],H^s)$ class. \\

We prove a slightly stronger version of Theorem~\ref{thm1} in
Section~\ref{locex}.

In the case $\alpha=1/2,$ we prove a result similar to
Theorem~\ref{subcritical}. However, the critical case is harder
and requires a new nonlocal maximum principle. We handle this case
in Section~\ref{glocrit}.

\begin{theorem}\label{thm2}
Assume that $\alpha=1/2,$ and that the initial data $u_0$ belongs
to $H^s$ with $s >1/2.$ Then there exists a global solution of
\eqref{bur1} which is real analytic in $x$ for any $t>0.$
\end{theorem}

The corresponding question for the quasi-geostrophic equation has
been a focus of significant effort (see e.g. \cite{CMT,CCW,CC})
and has been recently resolved independently and by different
means in \cite{KNV} and \cite{CafVas}. The proof in \cite{KNV} is
similar to the argument presented here.

Next, we prove the finite time blow up for the supercritical case
$\alpha<1/2.$

\begin{theorem}\label{thm3}
Assume that $0<\alpha<1/2.$ Then there exists smooth periodic
initial data $u_0(x)$ such that the solution $u(x,t)$ of
\eqref{bur1} blows up in $H^s$ for each $s>\frac32 -2\alpha$ in
finite time.
\end{theorem}
We will also obtain a fairly precise picture of blow up, similar
to that of Burgers equation without dissipation -- a shock is
formed where the derivative of the solution becomes infinite. In
the scenario we develop, the initial data leading to blow up is
odd and needs to satisfy a certain size condition, but no other
special assumptions. After this work has been completed, we became
aware of the preprint \cite{ADV}, where a result similar to
Theorem~\ref{thm3} is proved in the whole line (not periodic)
setting, and for a class of initial data satisfying certain
convexity assumption.

The blow up or global regularity for $\alpha <1/2$ remains open
for the surface quasi-geostrophic equation. The problem is that
the conservative surface quasi-geostrophic dynamics is not well
understood, in contrast to the non-viscous Burgers equation where
finite time shock formation is both well-known and simple.
Existence of blow up in the non-viscous quasi-geostrophic equation
remains a challenging open question (see e.g. \cite{Cord1}).
Still, some elements of the blow up construction we present here
may turn out to be useful in future attempts to attack the
question of blow up or regularity for the dissipative surface
quasi-geostrophic equation.

In Section~\ref{rough}, we prove existence of solutions with rough
initial data, $u_0 \in L^p,$ $1<p<\infty,$ when $\alpha=1/2$ (the
case $\alpha >1/2$ is similar). These solutions become smooth
immediately for $t>0,$ however the behavior near zero may be quite
singular. The uniqueness for such solutions is not known and
remains an interesting open problem.

The results of Theorems~\ref{subcritical}, \ref{thm1}, \ref{thm2},
\ref{thm3} remain valid for the case $s=3/2-2\alpha$, with slight
modifications. However, more subtle estimates are needed. This
critical space case is handled in Section~\ref{critical}.

\section{Local existence, uniqueness and regularity}\label{locex}

Denote by $P^N$ the orthogonal projection to the first $(2N+1)$
eigenfunctions of Laplacian, $e^{2\pi ikx},$ $k = 0, \pm 1, \dots,
\pm N.$ Consider Galerkin approximations $u^N(x,t),$ satisfying
\begin{equation}\label{gal}
u_t^N = P^N(u^N u^N_x) - (-\Delta)^\alpha u^N, \,\,\,u^N(x,0)=P^N
u_0(x).
\end{equation}
We start with deriving some a-priori bounds for the growth of
Sobolev norms.
\begin{lemma}\label{nlest}
Assume that $s\geq0$ and $\beta \geq 0.$ Then
\begin{equation}\label{nlestg}
\left| \int (u^N)^2 (-\Delta)^s u_x^N\,dx \right| \leq C\|u^N\|_q
\|u^N\|_{s+\beta}^2
\end{equation}
for any $q$ satisfying $q > 3/2 - 2\beta.$
\end{lemma}
\begin{proof}
On the Fourier side, the integral in \eqref{nlestg} is equal to
(up to a constant factor)
\[ \sum\limits_{k+a+b=0,|k|,|a|,|b| \leq N}
k|k|^{2s}\hat{u}^N(k)\hat{u}^N(a)\hat{u}^N(b)=:S. \]
Symmetrizing, we obtain
\begin{eqnarray}\label{sums11}
|S|=\frac13\left|\sum\limits_{k+a+b=0,|k|,|a|,|b| \leq N}
\left(k|k|^{2s}+a|a|^{2s}+b|b|^{2s}\right)\hat{u}^N(k)\hat{u}^N(a)\hat{u}^N(b)\right|\leq\\
\nonumber 2\sum\limits_{k+a+b=0,|a|\leq|b|\leq|k|\leq N}
\left|\,k|k|^{2s}+a|a|^{2s}+b|b|^{2s}\right||\hat{u}^N(k)\hat{u}^N(a)\hat{u}^N(b)|.
\end{eqnarray}
Next, note that under conditions $|a|\leq|b|\leq|k|,\ \ a+b+k=0$, we have
$|a|\leq|k|/2,\ |b|\geq|k|/2$ and
\begin{equation}\label{canc22}
 \begin{split} &
\left|k|k|^{2s}+a|a|^{2s}+b|b|^{2s}\right|=\left|b(|b|^{2s}-|b+a|^{2s})+a(|a|^{2s}-|k|^{2s})\right|\leq
\cr & C(s)|a||k|^{2s}\leq
C(s)|a|^{1-2\beta}|b|^{s+\beta}|k|^{s+\beta}.
\end{split}
\end{equation}
Thus
\begin{eqnarray*}
|S|\leq C\sum\limits_{k+a+b=0,|a|\leq|b|\leq|k|\leq N}
|a|^{1-2\beta}|b|^{s+\beta}|k|^{s+\beta}
|\hat{u}^N(k)\hat{u}^N(a)\hat{u}^N(b)|\leq\\
\nonumber C\|u^N\|_{s+\beta}^2\sum\limits_{|a|\leq N}
|a|^{1-2\beta}|\hat{u}^N(a)|\leq
C(\beta,q,s)\|u^N\|_q\|u^N\|_{s+\beta}^2.
\end{eqnarray*}
Here the second inequality is due to Parseval and convolution
estimate, and the third holds by H\"older's inequality for every
$q
> 3/2 - 2\beta.$
\end{proof}

Lemma~\ref{nlest} implies a differential inequality for the
Sobolev norms of solutions of \eqref{gal}.
\begin{lemma}\label{dest}
Assume that $\alpha>0,$ $q>3/2-2\alpha$, and $s\geq0$. Then
\begin{equation}\label{diffest22}
\frac{d}{dt} \|u^N\|_s^2 \leq C(q) \|u^N\|_q^{M(q,\alpha,s)}
-\|u^N\|^2_{s+\alpha}.
\end{equation}
If in addition $s=q$ then
\begin{equation}\label{diffest11}
\frac{d}{dt} \|u^N\|_s^2 \leq C(\epsilon)
\|u^N\|_s^{2+\frac{\alpha}{\epsilon}} -\|u^N\|^2_{s+\alpha},
\end{equation}
for any \begin{equation}\label{11epscon} 0<\epsilon< {\rm
min}\left(\frac{2q-3+4\alpha}{4},\alpha\right). \end{equation}
\end{lemma}
\begin{proof}
Multiplying both sides of \eqref{gal} by $(-\Delta)^s u^N,$ and
applying Lemma~\ref{nlest}, we obtain (here we put
$\beta:=\alpha-\epsilon,$ with $\epsilon$ satisfying
\eqref{11epscon})
\[
\frac{d}{dt} \|u^N\|_s^2 \leq C(q,\epsilon,\alpha,s) \|u^N\|_q
\|u^N\|_{s+\alpha-\epsilon}^2 -2\|u^N\|^2_{s+\alpha}.
\]
Observe that if $q \geq s+\alpha-\epsilon,$ the estimate
\eqref{diffest22} follows immediately. If $q<s+\alpha-\epsilon,$
by H\"older we obtain
\begin{equation}\label{hold11}
\|u^N\|^2_{s+\alpha-\epsilon}\leq\|u^N\|_{s+\alpha}^{2(1-\delta)}\|u^N\|_q^{2\delta}
\end{equation}
where $\delta=\dfrac{\epsilon}{s+\alpha-q}.$ Applying Young's
inequality we finish the proof of \eqref{diffest22} in this case.

The proof of \eqref{diffest11} is similar. We have
\[
\frac{d}{dt} \|u^N\|_s^2 \leq C(s,\epsilon,\alpha) \|u^N\|_s
\|u^N\|_{s+\alpha-\epsilon}^2 -2\|u^N\|^2_{s+\alpha}.
\]
Applying the estimate \eqref{hold11} with $q=s$
and $\delta=\epsilon/\alpha$ and Young's inequality we obtain
\[
\frac{d}{dt} \|u^N\|_s^2 \leq C\|u^N\|_s^{1+2\epsilon/\alpha}
\|u^N\|_{s+\alpha}^{2-2\epsilon/\alpha} -2\|u^N\|^2_{s+\alpha}\leq
C \|u^N\|_s^{2+\frac{\alpha}{\epsilon}} -\|u^N\|^2_{s+\alpha}.
\]
\end{proof}
The following lemma is an immediate consequence of
\eqref{diffest11} and local existence of the solution to the
differential equation $y'=Cy^{1+\alpha/2\epsilon},\ y(0)=y_0$.
\begin{lemma}\label{existenceini}
Assume $s>3/2-2\alpha$, $\alpha>0$ and $u_0 \in H^s.$ Then there
exists time $T=T(s,\alpha,\|u_0\|_s)$ such that for every $N$ we
have the bound (uniform in $N$)
\begin{equation}\label{cs2galini}
\|u^N\|_{s}(t)\leq C(s,\alpha,\|u_0\|_s),\ \ 0\leq t\leq T.
\end{equation}
\end{lemma}
\begin{proof}
From \eqref{diffest11}, we get that $z(t) \equiv \|u^N\|^2_s$
satisfies the differential inequality $z' \leq
Cz^{1+\alpha/2\epsilon}.$ This implies the bound \eqref{cs2galini}
for time $T$ which depends only on coefficients in the
differential inequality and initial data.
\end{proof}
Now, we obtain some uniform bounds for higher order $H^s$ norms of
the Galerkin approximations.
\begin{theorem}\label{existcinfgal}
Assume $s > 3/2-2\alpha$, $s\geq0$, $\alpha>0$ and $u_0 \in H^s.$
Then there exists time $T=T(s,\alpha,\|u_0\|_s)$ such that for
every $N$ we have the bounds (uniform in $N$)
\begin{equation}\label{cs2gal}
t^{n/2}\|u^N\|_{s+n\alpha}\leq C(n,s,\alpha,\|u_0\|_s),\ \ 0<
t\leq T,
\end{equation}
for any $n\geq0.$ Here time $T$ is the same as in
Lemma~\ref{existenceini}.
\end{theorem}
\begin{proof}
We are going to first verify \eqref{cs2gal} by induction for
positive integer $n$. For $n=0,$ the statement follows from
Lemma~\ref{existenceini}. Inductively, assume that $\|u^N\|^2_{s+n
\alpha}(t) \leq Ct^{-n}$ for $0 \leq t \leq T.$ Fix any
$t\in(0,T],$ and consider the interval $I=(t/2,t).$ By
\eqref{diffest22} with $s$ replaced by $s+n\alpha$ and $q$ by $s,$
we have for every $n \geq 0$
\begin{equation}\label{varest22}
\frac{d}{dt}\|u^N\|^2_{s+n\alpha} \leq C\|u^N\|_s^{M} -
\|u^N\|^2_{s+(n+1)\alpha}.
\end{equation}
Due to Lemma~\ref{existenceini} and our induction assumption,
\[ \int_{t/2}^t \|u^N\|^2_{s+(n+1)\alpha}\,ds \leq Ct +
C\|u^N(t/2)\|^2_{s+n\alpha} \leq Ct^{-n}. \] Thus we can find
$\tau \in I$ such that
\[ \|u^N(\tau)\|^2_{s+(n+1)\alpha} \leq C|I|^{-1}t^{-n} \leq
Ct^{-n-1}. \] Moreover, from \eqref{varest22} with $n$ changed to
$n+1$ we find that
\[ \|u^N(t)\|^2_{s+(n+1)\alpha} \leq
\|u^N(\tau)\|^2_{s+(n+1)\alpha} +Ct \leq Ct^{-n-1}, \] concluding
the proof for integer $n$.
 Non-integer $n$ can
be covered by interpolation:
\[\|u^N\|_{s+r\alpha}\leq\|u^N\|_s^{1-\frac{r}{n}}\|u^N\|_{s+n\alpha}^{\frac{r}{n}},\ \ \ 0<r\leq n.\]
\end{proof}

Now we are ready to prove existence and regularity of a weak
solution of Burgers equation \eqref{bur1}.
\begin{theorem}\label{existcinf}
Assume $s > 3/2-2\alpha$, $s\geq0$, $\alpha>0$, and $u_0 \in H^s.$
Then there exists $T(s,\alpha,\|u_0\|_s)>0$ and a solution
$u(x,t)$ of \eqref{bur1} such that
\begin{equation}\label{ws2}
u \in L^2([0,T], H^{s+\alpha}) \cap C([0,T],H^s);
\end{equation}
\begin{equation}\label{cs2}
t^{n/2}u  \in  C((0,T],H^{s+n \alpha})\cap L^\infty([0,T],H^{s+n
\alpha})
\end{equation}
for every $n>0.$
\end{theorem}

\begin{corollary}\label{cinfcor}
If $\alpha>0$ and $u_0 \in H^s$ with $s > 3/2 -2\alpha$, $s\geq0$,
then there exists a local solution $u(x,t)$ which is $C^\infty$
for any $0<t\leq T.$
\end{corollary}

\begin{proof}
The proof of Theorem~\ref{existcinf} is standard. It follows from
\eqref{gal} and \eqref{cs2gal} that for every small $\epsilon>0$
and every $r>0$ we have uniform in $N$ and $t\in[\epsilon,T]$
bounds
\begin{equation}\label{derb11}
\|u^N_t\|_r\leq C(r,\epsilon).
\end{equation}
By \eqref{cs2} and \eqref{derb11} and the well known compactness
criteria (see e.g. \cite{CF}, Chapter 8), we can find a
subsequence $u^{N_j}$ converging in $C([\epsilon,T], H^r)$ to some
function $u$. Since $\epsilon$ and $r$ are arbitrary one can apply
the standard subsequence of subsequence procedure to find a
subsequence (still denoted by $u^{N_j}$) which converges to $u$ in
$C((0,T], H^r)$, for any $r>0.$ The limiting function $u$ must
satisfy the estimates \eqref{cs2gal} and it is straightforward to
check that it solves the Burgers equation on $(0,T].$ Thus, it
remains to show that $u$ can be made to converge to $u_0$ strongly
in $H^s$ as $t\to0$.

We start by showing that $u$ converges to $u_0$ as $t \to 0$
weakly in $H^s.$ Let $\varphi(x)$ be arbitrary $C^\infty$
function. Consider
\[ g^N(t, \varphi) \equiv (u^N, \varphi) = \int u^N(x,t)\varphi(x)\,dx. \]
Clearly, $g^N( \cdot, \varphi) \in C([0,\tau]),$ where $\tau
\equiv T/2.$
Also, taking inner product of \eqref{gal} with $\varphi$ we can
show that for any $\delta >0,$
\begin{equation}\label{estim}
\int_0^\tau |g_t^N|^{1+\delta}\,dt \leq C\left(\int_0^\tau
\|u^N\|_{L^2}^{2+2\delta} \|\varphi\|^{1+\delta}_{W^1_\infty}\,dt
+ \int_0^\tau \|u^N\|^{1+\delta}_{L^2}
\|\varphi\|_{2\alpha}^{1+\delta}\,dt \right).
\end{equation}
Due to \eqref{diffest11}, the definition of $\tau,$ and the
condition $s \geq 0,$ we have that $\|u^N\|_{L^2} \leq C$ on
$[0,\tau],$ and thus $\|g^N_t(\cdot,\varphi)\|_{L^{1+\delta}} \leq
C(\varphi).$ Therefore the sequence $g^N(t,\varphi)$ is compact in
$C([0,\tau]),$ and we can pick a subsequence $g^{N_j}(t,\varphi)$
converging uniformly to a function $g(t,\varphi) \in C([0,\tau]).$
Clearly, by choosing an appropriate subsequence we can assume
$g(t,\varphi)=(u,\varphi)$ for $t \in (0,\tau].$ Next, we can
choose a subsequence $\{N_j\}$ such that $g^{N_j}(t,\varphi)$ has
a limit for any smooth function $\varphi$ from a countable dense
set in $H^{-s}.$ Given that we have uniform control over
$\|u^{N_j}\|_s$ on $[0,\tau],$ it follows that
$g^{N_j}(t,\varphi)$ converges uniformly on $[0,\tau]$ for every
$\varphi \in H^{-s}.$ Now for any $t>0,$
\begin{equation}\label{fin3w}
|(u-u_0,\varphi)|\leq
|(u-u^{N_j},\varphi)|+|(u^{N_j}-u^{N_j}_0,\varphi)|+|(u^{N_j}_0-u_0,\varphi)|.
\end{equation}
The first and the third terms in RHS of \eqref{fin3w} can be made
small uniformly in $(0,\tau]$ by choosing sufficiently large
$N_j$. The second term tends to zero as $t\to0$ for any fixed
$N_j$. Thus $u(\cdot,t)$ converges to $u_0(\cdot)$ as $t\to0$
weakly in $H^s$. Consequently,
\begin{equation}\label{upperhs}
\|u_0(\cdot)\|_{s}\leq \liminf_{t\to0}\|u(\cdot,t)\|_{s}.
\end{equation}
Furthermore, it follows from \eqref{diffest11} that for every $N$
the function $\|u^N\|_s^2(t)$ is always below the graph of the
solution of the equation
\[y_t=Cy^{1+\frac{\alpha}{2\epsilon}},\ \ y(0)=\|u_0\|_s^2.\]
By construction of the solution $u,$ the same is true for
$\|u\|_s^2(t)$. Thus,
$\|u_0\|_s\geq\limsup\limits_{t\to0}\|u\|_s(t)$. From this and
\eqref{upperhs}, we obtain that
$\|u_0\|_s=\lim\limits_{t\to0}\|u\|_s(t)$. This equality combined
with weak convergence finishes the proof.
\end{proof}


We next turn to the uniqueness. First of all, we obtain some
identities which hold for every weak solution of the Burgers
equation \eqref{bur1}. Let $u$ be a solution of the Burgers
equation in a sense of \eqref{burweak}, \eqref{burweak1}. Then for
any function $f$ of the form
\[f(x,t):=\sum\limits_{k=1}^K\varphi_k(x)\psi_k(t),\] where
$\varphi_k\in C^\infty({\mathbb T})$ and $\psi_k\in
C^\infty_0([0,T]),$ we have (see \eqref{burweak1})
\[
\sum\limits_{k=1}^K (u,\varphi_k)_t\psi_k=-\frac12(u^2,f_x)
-(u,(-\Delta)^{\alpha}f),\ \ \hbox{a.e.}\ t\in(0,T).
\]
Integrating and using integration by parts on the left hand side
we obtain
\begin{equation}\label{burweak2}
-\int\limits_0^T(u,f_t)\,dt=-\frac12\int\limits_0^T(u^2,f_x)\,dt
-\int\limits_0^T(u,(-\Delta)^{\alpha}f)\,dt.
\end{equation}
Applying closure arguments to \eqref{burweak2} and using inclusion
$u\in L^1([0,T],L^2)$ we derive the following statement.
\begin{lemma}\label{weak}
Let $u$ be a weak solution of the Burgers equation \eqref{bur1} in
the sense of \eqref{burweak}, \eqref{burweak1}. Then for every
function $f\in C^\infty_0([0,T],C^\infty({\mathbb{T}}))$ we have
\begin{equation}\label{burweak3}
-\int\limits_0^T(u,f_t)\,dt=-\frac12\int\limits_0^T(u^2,f_x)\,dt
-\int\limits_0^T(u,(-\Delta)^{\alpha}f)\,dt,\ \ \ f\in
C^\infty_0([0,T],C^\infty({\mathbb{T}})).
\end{equation}
\end{lemma}
Now, we are ready to prove
\begin{theorem}\label{unique}
Assume $v(x,t)$ is a weak solution of \eqref{bur1} for $0<\alpha
\leq 1/2,$ and initial data $u_0\in H^s$, $s>3/2-2\alpha$. If
\begin{equation}\label{smaluncon}
v(x,t) \in C([0,T],L^2)\cap L^{3/2\delta}([0,T], H^\delta), \ \ \
1\geq\delta>1/2,
\end{equation}
then $v(x,t)$ coincides with the solution $u(x,t)$ described in
Theorem~\ref{existcinf}.

\end{theorem}
\it Remark. \rm Theorems~\ref{unique} and \ref{existcinf} imply
Theorem~\ref{thm1}.
\begin{proof} We will need an auxiliary estimate for $\|v^2\|_{1-\delta}$. Recall that by integral
characterization of Sobolev spaces,
\[ \|v^2\|_{1-\delta} \leq C \left( \int_{\Sm^1} \int_{\Sm^1} \frac{|v(x)^2 -
v(y)^2|}{|x-y|^{1-\delta}}\,dxdy +\|v^2\|_{L^2} \right) \leq
C\|v\|_{L^\infty}\|v\|_{1-\delta}. \] Recall that for $1/2 <
\delta \leq 1,$ we have $\|v\|_{L^\infty} \leq C
\|v\|^{1-1/2\delta} \|v\|_\delta^{1/2\delta}.$ Applying this
inequality the $L^\infty$ norm and H\"older inequality to
$H^{1-\delta}$ norm above we obtain
\begin{equation}\label{raz}
\|v^2\|_{1-\delta}\leq C'\|v\|^{2-k}\|v\|_\delta^k,\ \ \
k:=\frac{3}{2\delta}-1.
\end{equation}
Now, let us obtain a bound for $\frac{dv}{dt}$. Since $v$ is a
weak solution of the equation \eqref{bur1}, due to
Lemma~\ref{weak} for every function $f(x,t)\in
C^\infty_0([0,T],C^\infty({\mathbb{T}}))$ we have
\begin{equation}\label{derivt}
\begin{split}
&
\left|\int_0^T(v,f_t)\,dt\right|\leq\frac12\int_0^T|(D^{1-\delta}(v^2),D^{\delta}f)|\,dt+\int_0^T
|(D^{2\alpha-\delta}v,D^{\delta}f)|\,dt\leq\cr &
\left(\int_0^T\|f\|_\delta^\gamma dt\right)^{1/\gamma}
\left(\frac12\left(\int_0^T\|v^2\|_{1-\delta}^{\gamma'}dt\right)^{1/\gamma'}+
\left(\int_0^T\|v\|_{1-\delta}^{\gamma'}dt\right)^{1/\gamma'}\right).
\end{split}
\end{equation}
Here $D:=(-\Delta)^{1/2}$, $\gamma:=\frac{3}{2\delta}$ and
$\gamma^{-1}+(\gamma')^{-1}=1$. It follows from \eqref{smaluncon},
\eqref{raz} and equality $\gamma'k=\gamma$ that the integral
$\int_0^T\|v^2\|_{1-\delta}^{\gamma'}dt$ is convergent. The
estimate for $\int_0^T\|v\|_{1-\delta}^{\gamma'}dt$ is similar and
even simpler. Thus, it follows from \eqref{derivt} that
$\frac{dv}{dt}$ belongs to $L^{\gamma'}([0,T],H^{-\delta})$.
Certainly, the same (and even more) is true for $\frac{du}{dt}$.
Thus, $((u-v)_t,(u-v))\in L^1([0,T])$. Moreover, it follows from
the definition of a weak solution, our assumptions and estimates
for $v_t$ and $v^2$ that for a.e. $t\in[0,T]$
\[
u_t -v_t = \frac12(u^2)_x-\frac12(v^2)_x-(-\Delta)^\alpha(u-v), \
\ \hbox{a.e.}\ t\in[0,T],
\]
where the equality is understood in $H^{-\delta}$ sense. Thus,
\begin{equation}\label{dva}
2((u-v)_t,(u-v))=((u^2)_x-(v^2)_x,u-v)- 2\|u-v\|_\alpha^2, \
\ \hbox{a.e.}\ t\in[0,T].
\end{equation}
For every fixed $t\in[0,T]$ where \eqref{dva} holds we approximate
$v$ in $H^{\delta}$ by smooth functions $v_n$. Direct calculations
using integration by parts give
\begin{eqnarray}\nonumber
&
((u^2)_x-(v_n^2)_x,u-v_n)=2(u_x,(u-v_n)^2)+2((u-v_n)_x,v_n(u-v_n))=\cr
&
2(u_x,(u-v_n)^2)-2((u-v_n)_x,(u-v_n)^2)+2((u-v_n)_x,u(u-v_n))=(u_x,(u-v_n)^2).
\end{eqnarray}
Due to \eqref{smaluncon} and \eqref{raz} we may apply closure
arguments in $H^\delta$ to obtain
\[((u^2)_x-(v^2)_x,u-v)=(u_x,(u-v)^2).\] Substituting into
\eqref{dva} we get
\begin{equation}\label{qqq111}
2((u-v)_t,(u-v))=(u_x,(u-v)^2)- 2\|u-v\|_\alpha^2, \ \
\hbox{a.e.}\ t\in[0,T].
\end{equation}

Note that $2((u-v)_t,(u-v))=\partial_t \|u-v\|^2$. Indeed, denote
$u-v=:g$. Recall that $g_t\in L^{\gamma'}([0,T],H^{-\delta})$ and
$g\in C([0,T],L^2)\cap L^{\gamma}([0,T], H^\delta)$. Approximate
$g$ in $W^1_{\gamma'}([0,T],H^{-\delta})\cap C([0,T],L^2)\cap
L^{\gamma}([0,T], H^\delta)$ by smooth functions $g_n$. Then
$2((g_n)_t,g_n)=\partial_t\|g_n\|^2$ and for every $t\in(0,T]$ we
have
\[\int\limits_0^t 2((g_n)_t,g_n)\,dt=\|g_n\|^2(t)-\|g_n\|^2(0).\] Now, we can take the limit to obtain
\[\int\limits_0^t 2(g_t,g)\,dt=\|g\|^2(t)-\|g\|^2(0).\] This
proves the desired identity.

Finally,
\[ |(u_x,(u-v)^2)|=\left| \int u_x(u-v)^2\,dx \right| \leq \|u-v\|^2_{L^{2p}}
\|u\|_{W^1_{p'}}, \] where $p=\frac{1}{1-\alpha}$ and $p'$ is the
H\"older conjugate exponent to $p$. The exponent $p$ was chosen so
that, by Gagliardo-Nirenberg inequality,
\[ \|u-v\|^2_{L^{2p}} \leq C\|u-v\|
\|u-v\|_\alpha.\] Also by Sobolev inequality
\[ \|u\|_{W^1_{p'}} \leq C\|u\|_r, \]
where $r =3/2-\alpha$.
Thus
\begin{equation}\label{qqq222}
 \left| \int u_x(u-v)^2\,dx
\right| \leq C\|u-v\|\|u-v\|_\alpha \|u\|_{r}.
\end{equation}
 Then from \eqref{qqq111}, \eqref{qqq222} we find
\begin{equation}\label{aeest}
\|u-v\|^2(t)-\|u-v\|^2(0)=2\int\limits_0^t((u-v)_s,(u-v))\,ds \leq
C\int\limits_0^t \|u-v\|^2 \|u\|^2_{r}\,ds,\ \
\hbox{for}\,\,\hbox{every} t\in[0,T].
\end{equation}
Now we are in position to apply Gronwall inequality.
Notice that $\int_0^T \|u\|_{s+\alpha}^2\,dt$ is controlled by
$\|u\|_{C([0,T],H^s)}^2$ because of \eqref{diffest11}. But if
$s>3/2-2\alpha,$ then $\int_0^T \|u\|^2_{r}\,dt$ is also under
control. This proves the Theorem.
\end{proof}


\section{Global existence and analyticity for the critical case
$\alpha=1/2$}\label{glocrit}

As we have already mentioned in the introduction, if $\alpha>1/2,$
then the local smooth solution can be extended globally. 
In Section~\ref{blowup}, we show that a blow up can happen in
finite time if $\alpha<1/2.$ Thus the only remaining case to
consider is the critical one, $\alpha=1/2.$ In this section, we
will prove

\begin{theorem}\label{critcase}
Assume $\alpha=1/2,$ and $u_0 \in H^s,$ $s>1/2.$ Then there exists
a global solution $u(x,t)$ of \eqref{bur1} which belongs to $C([0,\infty),H^s)\cap C((0,\infty),C^\infty)$.
If $v$ is another weak solution of \eqref{bur1}
with initial data $u_0$ such that
$v\in C([0,T],L^2)\cap L^{3/2\delta}([0,T],H^\delta)$ with some $\delta\in(1/2,1]$,
then $v$ coincides with $u$ on $[0,T]$.
\end{theorem}
We also state the following result separately to break the
otherwise unwieldy proof:
\begin{theorem}\label{analcrit}
The solution of Theorem~\ref{critcase} is real analytic for every
$t>0.$
\end{theorem}
Together, Theorems~\ref{critcase} and \ref{analcrit} imply
Theorem~\ref{thm2}.
We first discuss the global existence. Much of the discussion
follows \cite{KNV}; we reproduce the argument here for the sake of
completeness. Recall that a modulus of continuity is an arbitrary
increasing continuous concave function
$\omega\,:\,[0,+\infty)\to[0,+\infty)$ such that $\omega(0)=0$.
Also, we say that a function $f\,:\, \Rm \to \Rm$ has modulus of
continuity $\omega$ if $|f(x)-f(y)|\le\omega(|x-y|)$ for all
$x,y\in \Rm$.

The term $u u_x$ in the dissipative Burgers equation tends to make
the modulus of continuity of $u$ worse while the dissipation term
$(-\Delta)^{1/2} u$ tends to make it better. Our aim is to
construct some special moduli of continuity for which the
dissipation term always prevails and such that every periodic
$C^\infty$-function $u_0$ has one of these special moduli of
continuity.

Our moduli of continuity will be derived from one single function
$\omega(\xi)$ by scaling: $\omega_B(\xi)=\omega(B\xi).$ Note that
the critical ($\alpha=\frac12$) equation has a simple scaling
invariance: if $u(x,t)$ is a solution, then so is $u(Bx,Bt)$. This
means that if we prove that the modulus of continuity $\omega$ is
preserved by the evolution, then the whole family
$\omega_B(\xi)=\omega(B\xi)$ of moduli of continuity will also be
preserved (provided that we look at the initial data of all
periods). Also observe that if $\omega$ is unbounded, then every
$C^\infty$ periodic function has modulus of continuity $\omega_B$
if $B>0$ is sufficiently large.

We will eventually have an explicit expression for $\omega.$ For
now, we show how preservation of $\omega$ is used to control the
solution.
\begin{lemma}\label{dercon}
Assume that $\omega(\xi)$ satisfies
\begin{equation}\label{dercon1}
\omega'(0) < \infty, \,\,\,\omega''(0) = -\infty.
\end{equation}
Then if a smooth function $f$ has modulus of continuity $\omega,$
it must satisfy $\| f' \|_{L^\infty} < \omega'(0).$
\end{lemma}
\begin{proof}
Indeed, take a point $x\in\Rm$ at which $\max|f'|$ is attained and
consider the point $y=x+\xi$. Then we must have $f(y)-f(x)\le
\omega(\xi) $ for all $\xi\ge 0$. But the left hand side is at
least $|f'(x)|\xi-C\xi^2$ where $C=\frac 12\|f''\|_{L^\infty}$
while the right hand side can be represented as
$\omega'(0)\xi-\rho(\xi)\xi^2$ with $\rho(\xi)\to+\infty$ as
$\xi\to 0+$. Thus $|f'(x)|\le \omega'(0)- (\rho(\xi)-C)\xi$ for
all $\xi>0$ and it remains to choose some $\xi>0$ satisfying
$\rho(\xi)>C$.
\end{proof}

Given this observation, Theorem~\ref{critcase} will be proved as
follows. By Theorem~\ref{existcinf}, if $u_0 \in H^s,$ $s>1/2,$
then the solution immediately becomes $C^\infty$ and stays smooth
at least till time $T$. Hence it will preserve one of our moduli
of continuity, $\omega_B.$ This will imply that
$\|u_x\|_{L^\infty}$ remains bounded by $\omega'_B(0).$ Starting
at time $T$ we construct again Galerkin approximations which
define our solution for one more step in time. And we continue
this process inductively. Since our solution remains smooth and
its $W^{1,\infty}$ norm satisfies uniform bound, it implies
uniform boundedness of, for instance, $H^1$ norm for which our
local existence result is valid. Therefore our time step can be
chosen to be fixed. Thus we obtain smooth global solution of the
equation \eqref{bur1}.

We next proceed with the construction of $\omega$ and the proof of
preservation. Let us outline the only scenario how the modulus of
continuity satisfying \eqref{dercon1} can be lost.
\begin{lemma}\label{scen11}
Assume that a smooth solution $u(x,t)$ has modulus of continuity
$\omega$ at some time $t_0.$ The only way this modulus of
continuity may be violated is if there exists $t_1 \geq t_0$ and
$y,z,\,y\not=z$, such that $u(y,t_1)-u(z,t_1)=\omega(|y-z|),$
while for all $t < t_1,$ the solution has modulus of continuity
$\omega.$
\end{lemma}
\begin{proof}
Assume that $u(x,t)$ loses modulus of continuity $\omega.$ Define
\[ \tau = {\rm sup}\{ t: \, \forall x,y, \,\,\,\,|u(x,t)-u(y,t)| \leq
\omega(|x-y|). \} \] Then $u$ remains smooth up to $\tau,$ and, by
local existence and regularity theorem, for a short time beyond
$\tau.$ Suppose that $|u(x,\tau)-u(y,\tau)| < \omega(|x-y|)$ for
all $x \ne y.$ We claim that in this case $u$ has modulus of
continuity $\omega$ for all $t>\tau$ sufficiently close to $\tau.$
Indeed, by Lemma~\ref{dercon} at the moment $\tau$ we have
$\|u'\|_{L^\infty} < \omega'(0).$ By continuity of derivatives and
compactness in space variable, this also holds for $t>\tau$ close
to $\tau,$ which immediately takes care of the inequality
$|u(x,t)-u(y,t)| < \omega(|x-y|)$ for small $|x-y|.$ Observe that
we only need to consider $x,y$ within a fixed bounded domain since
$u$ is periodic and $\omega$ increasing. Thus, it suffices to show
that $|u(x,t)-u(y,t)| < \omega(|x-y|)$ holds for all $t$ close
enough to $\tau$ and $x,y$ such that $\delta \leq |x-y| \leq
\delta^{-1}$ with some $\delta >0.$ But this follows immediately
from the inequality for time $\tau,$ smoothness of the solution
and compactness of the domain.
\end{proof}
\it Remark. \rm The key point of the above lemma is that we do not
have to worry about $\omega$ being violated "first" at the
diagonal $x=y,$ namely $||u'||_{L^{\infty}} = \omega'(0);$ the
modulus of continuity equality must happen first at two distinct
points. This knowledge makes the argument below simpler by ruling
out the extra case which otherwise would have to be considered.

\begin{proof}[Proof of Theorem~\ref{critcase}]
Assume now that $u(y,t_1)-u(z,t_1)=\omega(|y-z|)$ for some $y,z$
and $|y-z|=\xi>0.$ We will henceforth omit $t_1$ from notation.
The plan now is to show that we have necessarily
$\frac{d}{dt}(u(y,t_1)-u(z,t_1))<0.$ 
We need to estimate the flow and the dissipative terms entering
the Burgers equation. First, note that
\[ u(y)u'(y) = \left. \frac{d}{dh} u(y+hu(y))\right|_{h=0} \]
and similarly for $z.$ But
\[ u(y+hu(y)) - u(z+hu(z)) \leq \omega\left(|y-z|+h|u(y)-u(z)|\right)\leq
\omega(\xi +h\omega(\xi)). \] Since also $u(y)-u(z)=\omega(\xi),$
we conclude that
\[ u(y)u'(y)-u(z)u'(z) \leq \omega(\xi)\omega'(\xi). \]
Note that we assume differentiability of $\omega$ here. The
$\omega$ that we will construct below is differentiable except at
one point, and this special point is handled easily (by using the
larger of the one-sided derivatives). Next let us estimate the
difference of dissipative terms. Due to translation invariance, it
is sufficient to consider $y=\xi/2$ and $z=-\xi/2.$ Let us denote
by $P_h$ the one dimensional Poisson kernel, $P_h(x) = \dfrac1\pi
\dfrac{h}{x^2+h^2}.$ Recall that
\[ -(-\Delta)^{1/2}u(x) = \left. \frac{d}{dh} P_h * u(x) \right|_{h=0}. \]
By the Poisson summation formula, this equality is valid for
periodic $u(x)$ of every period. By symmetry and monotonicity of
the Poisson kernel,
\begin{eqnarray*} (P_h*u)(y)-(P_h*u)(z) =
\int_0^\infty
[P_h(\xi/2-\eta)-P_h(-\xi/2-\eta)](u(\eta)-u(-\eta))\,d\eta \leq
\\ \int_0^\infty
[P_h(\xi/2-\eta)-P_h(-\xi/2-\eta)]\omega(2\eta)\,d\eta = \\
\int_0^\xi P_h(\xi/2-\eta)\omega(2\eta)\,d\eta + \int_0^\infty
P_h(\xi/2+\eta)[\omega(2\eta+2\xi)-\omega(2\eta)]\,d\eta.
\end{eqnarray*}
The last formula can also be rewritten as
\[ \int_0^{\xi/2} P_h(\eta)[\omega(\xi+2\eta) +
\omega(\xi-2\eta)]\,d\eta + \int_{\xi/2}^\infty P_h(\eta)
[\omega(2\eta+\xi)-\omega(2\eta-\xi)]\,d\eta. \] Since
$\int_0^\infty P_h(\eta)\,d\eta =1/2,$ we see that the difference
$(P_h*u)(y)-(P_h*u)(z) -\omega(\xi)$ can be estimated from above
by
\[ \int_0^{\xi/2} P_h(\eta)[\omega(\xi+2\eta) +
\omega(\xi-2\eta)-2\omega(\xi)]\,d\eta + \int_{\xi/2}^\infty
P_h(\eta)
[\omega(2\eta+\xi)-\omega(2\eta-\xi)-2\omega(\xi)]\,d\eta.
\]
Dividing by $h$ and passing to the limit as $h \rightarrow 0+,$ we
obtain the following upper bound on the contribution of the
dissipative term into the time derivative
\[ \frac1\pi\int_0^{\xi/2} \frac{\omega(\xi+2\eta) +
\omega(\xi-2\eta)-2\omega(\xi)}{\eta^2}\,d\eta + \frac1\pi
\int_{\xi/2}^\infty
\frac{\omega(2\eta+\xi)-\omega(2\eta-\xi)-2\omega(\xi)}{\eta^2}\,d\eta.
\]
Note that due to concavity of $\omega,$ both terms are strictly
negative. We will denote the first integral by $I_{\omega,1}(\xi)$
and the second integral by $I_{\omega,2}(\xi).$

Now we are ready to define our modulus of continuity. Let us set
$\xi_0 \equiv \left(\frac{K}{4\pi}\right)^2,$ where $K$ is to be
chosen later. Then $\omega$ is given by
\begin{equation}\label{modcon}
\omega(\xi) = \left\{ \begin{array}{ll} \frac{\xi}{1+K\sqrt{\xi}}
& {\rm for}\,\,\,0 \leq \xi \leq
\xi_0;\\
C_K \log \xi & {\rm for}\,\,\,\xi \geq \xi_0. \end{array} \right.
\end{equation}
Here $C_K$ is chosen to provide continuity of $\omega.$  A direct
computation shows that
\begin{equation}\label{Ck}
C_K \sim (\log K)^{-1}\,\,\,{\rm as}\,\,\, K \rightarrow \infty.
\end{equation}
One can check that if $K$ is sufficiently large, then $\omega$ is
concave, with negative and increasing second order derivative on
both intervals in \eqref{modcon} (on the first interval,
$\omega''(\xi)= -K(3\xi^{-1/2}+K)/4(1+K\sqrt{\xi})^3$). The first
derivative of $\omega$ may jump at $\xi_0,$ but the left
derivative at $\xi_0$ is $\sim K^{-2},$ while the right derivative
is $\sim K^{-2} (\log K)^{-1}.$ We choose $K$ large enough so that
the left derivative is larger than the right derivative assuring
concavity. The bound on the flow term is valid at $\xi_0$ using
the value of larger (left) derivative.

It remains to verify that we have
\[ \omega(\xi)\omega'(\xi) + I_{1,\omega}(\xi) +
I_{2,\omega}(\xi) \leq 0 \] for any $\xi.$

I. \it The case $\xi \leq \xi_0.$ \rm Using the second order
Taylor formula and the fact that $\omega''$ is negative and
monotone increasing on $[0,\xi],$ we obtain that
\[ \omega(\xi+2\eta)+\omega(\xi-2\eta) \leq
\omega(\xi)+\omega'(\xi)2\eta+\omega(\xi)-\omega'(\xi)2\eta+2\omega''(\xi)\eta^2.
\]
This leads to an estimate
\[ I_{1,\omega}(\xi) \leq \frac{1}{\pi} \xi \omega''(\xi). \]
From \eqref{modcon}, we find that
\[ 2\omega(\xi)\omega'(\xi) = \frac{\xi^{1/2}(2\xi^{1/2}
+K\xi)}{(1+K\xi^{1/2})^3}, \] while
\[ \frac1\pi \xi \omega''(\xi) =
-\frac{K(3\xi^{1/2}+K\xi)}{4\pi(1+K\xi^{1/2})^3}. \] Taking into
account that $\xi \leq \xi_0,$ we find
\begin{equation}\label{for1}
2\omega(\xi)\omega'(\xi) +I_{1,\omega}(\xi) \leq0,
\end{equation}
for any $K$.

II. \it The case $\xi > \xi_0.$ \rm 
Due to concavity, we have
$\omega(\xi+2\eta) \leq \omega(2\eta-\xi)+\omega(2\xi),$ and thus
\[ I_{2,\omega}(\xi) \leq \frac1\pi \int_{\xi/2}^\infty
\frac{\omega(2\xi)-2\omega(\xi)}{\eta^2}\,d\eta. \] Clearly we
have $\omega(2\xi)\leq \frac32 \omega(\xi)$ for $\xi \geq \xi_0$
provided that $K$ was chosen large enough. In this case, we obtain
$I_{2,\omega}(\xi) \leq -\frac{\omega(\xi)}{\pi \xi}.$ Now it
follows from \eqref{modcon} that $2\omega(\xi)\omega'(\xi) =
2C_K^2\xi^{-1}\log \xi,$ while $\xi^{-1}\omega(\xi) =
C_K\xi^{-1}\log \xi.$ Given \eqref{Ck}, it as clear that
\begin{equation}\label{for2}
2\omega(\xi)\omega'(\xi)+I_{2,\omega}(\xi) \leq0,\ \ \ \xi\geq
\xi_0
\end{equation}
if only $K$ was chosen sufficiently large.
\end{proof}
Observe that as a byproduct, the proof also yields uniform in time
control of $\|u'\|_{L^\infty}.$
\begin{corollary}\label{quantcon}
Assume that the initial data $u_0(x)$ is such that $\|
u_0'\|_{L^\infty} < \infty.$ Then for every time $t,$ the solution
$u(x,t)$ of the critical Burgers equation satisfies
\[ \| u'(x,t) \|_{L^\infty} \leq \| u_0'\|_{L^\infty}
\exp(C\|u_0\|_{L^\infty}). \]
\end{corollary}
\begin{proof}
Choose $B$ so that $u_0(x)$ has the modulus of continuity
$\omega_B.$ Given the asymptotic behavior of $\omega$ for large
$\xi,$ this is guaranteed if
\[ C_K \log \left( \frac{B}{\|u_0'\|_{L^\infty}} \right) \geq
\|u_0\|_{L^\infty}. \]
The Corollary then follows from
\eqref{modcon} and preservation of $\omega_B$ by evolution.
\end{proof}

Finally, we prove Theorem~\ref{analcrit}, establishing analyticity
of the solution.
\begin{proof}[Proof of Theorem~\ref{analcrit}]
We will assume that the initial data $u_0\in H^{2}$. Even if we
started from $u_0 \in H^s,$ $s >1/2,$ Theorem~\ref{existcinf}
implies that we gain the desired smoothness immediately.

Let us rewrite the equation \eqref{gal} on the Fourier side
($\alpha=1/2$, without loss of generality assume the period is
equal to one):
\[ \frac{d\hat{u}^N(k)}{dt}=\pi i\sum\limits_{a+b=k,\,|a|,|b|,|k|\leq N}k\hat{u}^N(a) \hat{u}^N(b)
- |k|\hat{u}^N(k). \] To simplify notation we will henceforth omit
the restrictions $|a|,|b|,|k|\leq N$ in the summation, but they
are always present in the remainder of the proof. Put
$\xi_k^N(t):=\hat{u}^N(k,t)e^{\frac12|k|t}$. Observe that since
$u(x,t)$ is real, $\overline{\xi}^N_k = \xi^N_{-k}.$ We have
\begin{equation}\label{foureq}
\frac{d\xi_k^N}{dt}=\pi
i\sum\limits_{a+b=k}e^{-\gamma_{a,b,k}t}k\xi_a^N \xi_b^N -
\frac12|k|\xi_k^N,
\end{equation}
where $\gamma_{a,b,k}:=\frac12(|a|+|b|-|k|)$. Note that
\begin{equation}\label{gamma}
0\leq\gamma_{a,b,k}\leq\min\{|a|,\,|b|\}.
\end{equation}

Consider $Y_N(t):=\sum\limits_k|k|^4|\xi_k^N(t)|^2$. Then we have
\begin{equation}\label{Y}
\begin{split}
& \frac{dY_N}{dt}=\Re\left({-2\pi
i\sum\limits_{a+b+k=0}e^{-\gamma_{a,b,k}t}k|k|^4\xi_a^N
\xi_b^N\xi_k^N}\right) - \sum\limits_k|k|^5|\xi_k^N|^2 \cr &
=\Re\left({-2\pi i\sum\limits_{a+b+k=0}k|k|^4\xi_a^N
\xi_b^N\xi_k^N}\right)+\Re\left({-2\pi
i\sum\limits_{a+b+k=0}(e^{-\gamma_{a,b,k}t}-1)k|k|^4\xi_a^N
\xi_b^N\xi_k^N}\right) \cr & -
\sum\limits_k|k|^5|\xi_k^N|^2=:I_1+I_2+I_3.
\end{split}
\end{equation}
Symmetrizing $I_1$ over $a,\ b$ and $k$ we obtain
\begin{equation}
I_1=\frac{2\pi}{3}\Re\left({-i\sum\limits_{a+b+k=0}(k|k|^4+a|a|^4+b|b|^4)\xi_a^N
\xi_b^N\xi_k^N}\right).
\end{equation}
Thus
\begin{eqnarray}\label{I_1}
 |I_1|\leq 4\pi \sum\limits_{a+b+k=0,\,\,
|a|\leq|b|\leq|k|}|k|k|^4+a|a|^4+b|b|^4||\xi_a^N||\xi_b^N||\xi_k^N|\leq
\\ \nonumber \leq 160\pi \sum\limits_{a+b+k=0,\,\,
|a|\leq|b|\leq|k|}|a||b|^2|k|^2|\xi_a^N||\xi_b^N||\xi_k^N|\\
\nonumber\leq 160 \pi Y_N\sum|a||\xi_a^N|\leq C_1 Y_N^{3/2}.
\end{eqnarray}
Here in the second step we used $a+b+k=0$ (compare to
\eqref{canc22}), and in the last step we used H\"older inequality:
\begin{equation}\label{hold222}
 \sum\limits_a |a||\xi_a^N| \leq \left(\sum\limits_{a \ne 0}
|a|^{-2} \right)^{1/2} Y_N^{1/2}(t).
\end{equation}
 For $I_2$ we have
\[ |I_2|\leq 2\pi \sum\limits_{a+b+k=0}{\rm
 min}(|a|,|b|)t|k|^5|\xi_a^N||\xi_b^N||\xi_k^N|. \]
Here we used \eqref{gamma}. Furthermore,
\begin{eqnarray*}
\sum\limits_{a+b+k=0}{\rm
 min}(|a|,|b|)|k|^5|\xi_a^N||\xi_b^N||\xi_k^N| \leq
\sum\limits_{a+b+k=0,\,\, |a|\leq|b|\leq|k|}3|a||k|^5
|\xi_a^N||\xi_b^N||\xi_k^N|+ \\ \sum\limits_{a+b+k=0,\,\,
|b|\leq|a|\leq|k|}3|b||k|^5 |\xi_a^N||\xi_b^N||\xi_k^N| \leq
\sum\limits_{a+b+k=0}6|a||b|^{5/2}|k|^{5/2}|\xi_a^N||\xi_b^N||\xi_k^N|
\\ \leq 6 \left(\sum\limits_a |a||\xi_a^N| \right) \left(\sum\limits_k |k|^5
|\xi^N_k|^2 \right).
\end{eqnarray*}
We used Young's inequality for convolution in the last step.
Combining all estimates and applying \eqref{hold222}, we obtain
\begin{equation}\label{I_2}
 |I_2|\leq CtY_N^{1/2}\sum\limits_k |k|^5 |\xi^N_k|^2.
\end{equation}
 Combining \eqref{Y}, \eqref{I_1} and
\eqref{I_2} we arrive at
\begin{equation}\label{estY}
\frac{dY_N}{dt}\leq C_1 Y_N^{3/2}+(C_2
Y_N^{1/2}t-1)\sum\limits_k|k|^5|\xi_k^N|^2.
\end{equation}

Note that $Y_N(0)=\|u^N_0\|^2_2$. Thus we have a differential
inequality for $Y_N$ ensuring upper
bound on $Y_N$ uniform in $N$ 
for a short time interval $\tau$ which depends only on
$\|u_0\|_2$. Observe that Theorem~\ref{existcinfgal} and
Corollary~\ref{quantcon} ensure that the $H^2$ norm of solution
$u(x,t)$ is bounded uniformly on $[0,\infty).$ Thus we can use the
above construction to prove for every $t_0>0$ uniform in $N$ and
$t>t_0$ bound on $\sum_k |\hat{u}^N(k,t)|^2 e^{\delta |k|}$ for
some small $\delta(t_0,u_0)>0.$ By construction of $u$, it must
satisfy the same bound.
\end{proof}


\section{Blow-up for the supercritical case.}\label{blowup}
\medskip

Our main goal here is to prove Theorem~\ref{thm3}. At first, we
are going to produce smooth initial data $u_0(x)$ which leads to
blow up in finite time in the case where the period $2L$ is large.
After that, we will sketch a simple rescaling argument which gives
the blow up for any (and in particular unit) period.

The proof will be by contradiction. We will fix $L$ and the
initial data, and assume that by time $T=T(\alpha)$ the blow up
does not happen. In particular, this implies that there exists $N$
such that $\|u(x,t)\|_{C^3} \leq N$ for $0 \leq t \leq T.$ This
will lead to a contradiction. The overall plan of the proof is to
reduce the blow up question for front-like data to the study of a
system of differential equations on the properly measured
steepness and size of the solution. To control the solution, the
first tool we need is a time splitting approximation. Namely,
consider a time step $h,$ and let $w(x,t)$ solve
\begin{equation}\label{weq}
w_t = ww_x, \,\,\,w(x,0)=u_0(x),
\end{equation}
while $v(x,t)$ solves
\begin{equation}\label{veq}
v_t = -(-\Delta)^{\alpha}v, \,\,\,v(x,0)=w(x,h).
\end{equation}
The idea of approximating $u(x,t)$ with time splitting is fairly
common and goes back to the Trotter formula in the linear case
(see for example \cite{BM}, page 120, and \cite{T3}, page 307, for
some applications of time splitting in nonlinear setting). The
situation in our case is not completely standard, since the
Burgers equation generally does blow up, and moreover the control
we require is in a rather strong norm.

The solution of the problem \eqref{veq} with the initial data $v_0(x)$
is given by the convolution
\begin{equation}\label{disssol} v(x,t)
= \int_\Rm \Phi_t(x-y) v_0(y)\,dy\ (=e^{-(-\Delta)^\alpha t}v_0(x)),
\end{equation}
where
\begin{equation}\label{Phi1}
\Phi_t(x) = t^{-1/2\alpha} \Phi(t^{-1/2\alpha}x), \,\,\,\Phi(x) =
\frac{1}{2\pi}\int_\Rm \exp (ix\xi-|\xi|^{2\alpha}) \,d\xi.
\end{equation}
It is evident that $\Phi(x)$ is even and $\int \Phi(x)\,dx=1.$  We
will need the following further properties of the function $\Phi:$
\begin{equation}\label{Phi}
\Phi(x) >0;\,\,\,x\Phi'(x) \leq 0,\,\,\,\Phi(x) \leq
\frac{K(\alpha)}{1+|x|^{1+2\alpha}},\,\,\,\left|\Phi'(x)\right|
\leq \frac{K(\alpha)}{1+|x|^{2+2\alpha}}.
\end{equation}
These properties are not difficult to prove; see e.g.
\cite{Feller} for some results, in particular positivity (Theorem
XIII.6.1). 
We need the following lemma.

\begin{lemma}
For every $f\in C^{n+1},\ \ n\geq0$,
\begin{equation}\label{expcon}
\|(e^{-(-\Delta)^\alpha t} -1)f\|_{C^n} \leq C(\alpha) t
\|f\|_{C^{n+1}}.
\end{equation}
\end{lemma}
\begin{proof}
Obviously, it is sufficient to prove the Lemma for $n=0$. We have (see \eqref{disssol}, \eqref{Phi1}, \eqref{Phi})
\begin{eqnarray}
& \left|(e^{-(-\Delta)^\alpha t}
-1)f\right|=\left|\int\limits_{-\infty}^{\infty}\Phi_t(y)(f(x-y)-f(x))dy\right|\leq\cr &
\left|\int\limits_{-1}^{1}\Phi_t(y)(f(x-y)-f(x))dy\right|+
\left|\int\limits_{|y|\geq1}\Phi_t(y)(f(x-y)-f(x))dy\right|\leq\cr &
2\|f\|_{C^1}t^{\frac{1}{2\alpha}}\int\limits_0^{t^{-\frac{1}{2\alpha}}}y\Phi(y)dy+
4\|f\|_{C}\int\limits_{t^{-\frac{1}{2\alpha}}}^\infty\Phi(y)dy\leq
C\frac{K(\alpha)}{\alpha}\|f\|_{C^1}t.
\end{eqnarray}
\end{proof}
Next lemma provides local solvability for our splitting system.
\begin{lemma}\label{localforh}
Assume $\|u_0(x)\|_{C^3} \leq N.$ Then for all $h$ small enough,
$v(x,h)$ is $C^3$ and is uniquely defined by \eqref{weq},
\eqref{veq}. Moreover, it suffices to assume $h \leq CN^{-1}$ to
ensure
\begin{equation}\label{wvcon}
\|w(x,t)\|_{C^3}, \|v(x,t)\|_{C^3} \leq 2N
\end{equation}
for $0 \leq t \leq h.$
\end{lemma}
\begin{proof}
Using the characteristics one can explicitly solve equation
\eqref{weq}. We have $w(t,y)=u_0(x)$, where $x=x(y)$ is such that
\begin{equation}
y=x-u_0(x)t.
\end{equation}
Now, implicit function theorem and direct computations show that
$\|w(t,\cdot)\|_{C^3}\leq 2\|u_0\|_{C^3}$ provided that
$\|u_0\|_{C^3}t\leq c$ for some small constant $c>0$. This proves
the statement of the Lemma for $w$. To prove it for $v$ we just
notice that $v$ is a convolution of the $w(h,x)$ with $\Phi_t(x)$.
Since $\|\Phi_t\|_{L^1}=1$ we obtain that
$\|v(t,\cdot)\|_{C^3}\leq\|w(h,\cdot)\|_{C^3}.$
\end{proof}

The main time
splitting result we require is the following
\begin{proposition}\label{timesplit}
Assume that 
$\|u_0(x)\|_{C^3} \leq N$ for $0 \leq t \leq T.$ Define $v(x,t)$
by \eqref{weq}, \eqref{veq} with time step $h.$ Then for all $h$
small enough, we have
\[ \|u(x,h)-v(x,h)\|_{C^1} \leq
C(\alpha,N)h^2.\]
\end{proposition}
\begin{proof}
Since $\|u_0(x)\|_{C^3} \leq N,$ let us choose $h$ as in
Lemma~\ref{localforh}. Notice that by Duhamel's principle,
\[ u(x,h) = e^{-(-\Delta)^\alpha h} u_0(x) + \int\limits_0^h
e^{-(-\Delta)^\alpha (h-s)} (u(x,s) u_x(x,s))\,ds, \] while
\[ v(x,h) = e^{-(-\Delta)^\alpha h} u_0(x) + \int\limits_0^h
e^{-(-\Delta)^\alpha h} (w(x,s) w_x(x,s))\,ds. \]
Then it follows from \eqref{expcon} that
\begin{eqnarray}\nonumber
\|u(x,h)-v(x,h)\|_{C^1} \leq \int\limits_0^h \| e^{-(-\Delta)^\alpha (h-s)}
(u(x,s) u_x(x,s)) - e^{-(-\Delta)^\alpha h} (w(x,s)
w_x(x,s))\|_{C^1}\,ds \leq \\
\nonumber \int\limits_0^h \|u(x,s) u_x(x,s) - w(x,s)
w_x(x,s)\|_{C^1} \,ds + \int\limits_0^h
\|\left(e^{-(-\Delta)^\alpha (h-s)}-1\right) u(x,s)
u_x(x,s)\|_{C^1} + \\ \nonumber \int\limits_0^h
\|\left(e^{-(-\Delta)^\alpha h}-1\right) w(x,s)
w_x(x,s)\|_{C^1}\,ds \leq \int\limits_0^h \|u(x,s) u_x(x,s) -
w(x,s) w_x(x,s)\|_{C^1}\,ds +
\\ \label{uvdiff} C(\alpha)h \int\limits_0^h \left( \|u(x,s)u_x(x,s)\|_{C^2}+
\|w(x,s) w_x(x,s)\|_{C^2} \right)\,ds.
\end{eqnarray}
From \eqref{wvcon}, it follows that the last integral does not
exceed $C(\alpha) N^2 h^2.$ To estimate the remaining integral, we
need the following
\begin{lemma}\label{uvrough}
For every $0 \leq s \leq h,$ we have $\|u(x,s)-w(x,s)\|_{C^2} \leq
C(\alpha)N^2 h.$
\end{lemma}
\begin{proof}
Observe that $g(x,s) \equiv u(x,s) - w(x,s)$ solves
\[ g_t = gu_x +wg_x -(-\Delta)^\alpha u, \,\,\,g(x,0)=0. \]
Thus
\[ g(x,t) = \int\limits_0^t (gu_x +wg_x -(-\Delta)^\alpha u)\,ds.
\]
Because of \eqref{wvcon} and the assumption on $u,$ we have
$\|gu_x\|_{C^2},\|wg_x\|_{C^2} \leq CN^2,$ and $\|(-\Delta)^\alpha
u\|_{C^2} \leq CN.$ 
Therefore, we can estimate that $\|g(x,t)\|_{C^2} \leq
C(\alpha)tN^2,$ for every $0 \leq t \leq h.$
\end{proof}
From Lemma~\ref{uvrough} it follows that
\begin{eqnarray*}
\int\limits_0^h \|u u_x - w w_x\|_{C^1}\,ds \leq \int\limits_0^h
\left(\|(u-w) u_x\|_{C^1}+\|w(u_x-w_x)\|_{C^1} \right)\,ds \leq
C(\alpha)N^3 h^2.
\end{eqnarray*}
This completes the proof of Proposition~\ref{timesplit}.
\end{proof}

The next stage is to investigate carefully a single time splitting
step. The initial data $u_0(x)$ will be smooth, $2L$-periodic,
odd, and satisfy $u_0(L)=0.$ It is not hard to see that all these
assumptions are preserved by the evolution. We will assume a
certain lower bound on $u_0(x)$ for $0 \leq x \leq L,$ and derive
a lower bound that must hold after the small time step. The lower
bound will be given by the following piecewise linear functions on
$[0,L]:$
\[ \varphi (\kappa,H,a,x) = \left\{ \begin{array}{ll} \kappa x,
&  0 \leq x \leq \delta \equiv H/\kappa \\
H, & \delta \leq x \leq L-a \\
\frac{H}{a} (L-x), & L-a \leq x \leq L. \end{array} \right.
\]

Here $L,$ $\kappa,$ $H$ and $a$ may depend only on $\alpha$ and
will be specified later. We will set $a\leq L/4,$ $\delta\leq L/4$
and will later verify that this condition is preserved throughout
the construction. We assume that blow up does not happen until
time $T$ (to be determined later). Let $N = {\rm sup}_t
\|u(x,t)\|_{C^3}.$
\begin{lemma}\label{burstep}
Assume that the initial data $u_0(x)$ for the equation \eqref{weq}
satisfies the above assumptions. Then for every $h$ small enough
($h \leq CN^{-1}$ is sufficient), we have
\[ w(x,h) \geq \varphi \left(\frac{\kappa}{1-\kappa h}, H, a+\|u_0\|_{L^\infty}h,x\right),\ \ \ 0\leq x\leq L. \]
\end{lemma}
\begin{proof}
The Burgers equation can be solved explicitly using
characteristics. The existence of $C^3$ solution $w(x,t)$ for $t
\leq h$ is assured by the assumption on the initial data and $h.$
\end{proof}

Now we consider the effect of the viscosity time step. Suppose that the initial data $v_0(x)$ for \eqref{veq}
satisfies the same conditions as stated for $u_0(x)$ above:
periodic, odd,
$v_0(L)=0.$ 
Then we have
\begin{lemma}\label{disstep}
Assume that for $0 \leq x \leq L,$ $v_0(x) \geq \varphi(\kappa, H,
a, x).$ Moreover, assume that
\begin{equation}\label{condstep}
H\kappa^{-1}\leq a,\ \ \ L\geq 4a,\ \ \
L^{-2\alpha}\|v_0\|_{L^\infty} \leq 4Ha^{-2\alpha}.
\end{equation}
Then for every sufficiently small $h,$ we have
\[ v(x,h) \geq \varphi(\kappa(1-C(\alpha)h H^{-2\alpha}\kappa^{2\alpha}),
H(1-C(\alpha)hH^{-2\alpha}\kappa^{2\alpha}), a,x),\ \ \ 0\leq x\leq L.
\]
\end{lemma}
\begin{proof}
Let us compute
\begin{eqnarray}\nonumber
 v(x,h) = \int_{-\delta}^\delta \Phi_h(x-y)
v_0(y)\,dy+ \int_\delta^{L-a}
(\Phi_h(y-x)-\Phi_h(x+y))v_0(y)\,dy+ \\
\label{smallx} \int_{L-a}^\infty
(\Phi_h(y-x)-\Phi_h(x+y))v_0(y)\,dy.
\end{eqnarray}
In the last integral in \eqref{smallx}, we estimate by Mean Value
Theorem $|\Phi_h(y-x)-\Phi_h(y+x)| \leq 2x
\left|\Phi_h'(\tilde{y})\right|,$ where $\tilde{y} \in (y-x,y+x).$
Using \eqref{Phi}, we see that the last integral in \eqref{smallx}
is controlled by $C(\alpha)hxL^{-1-2\alpha}\|v_0\|_{L^\infty}.$
The second integral on the right hand side of \eqref{smallx} can
be estimated from below by
\begin{equation}\label{plato} H
\int_\delta^\infty (\Phi_h(y-x)-\Phi_h(x+y))\,dy
\end{equation}
with an error, which, by the previous computation, does not exceed
$C(\alpha)hxL^{-1-2\alpha}\|v_0\|_{L^\infty}.$ The expression in
\eqref{plato} is equal to
\begin{equation}\label{plato2}
H \int_{\delta -x}^{\delta+x} \Phi_h(z)\,dz.
\end{equation}
For the first integral in \eqref{smallx} we have
\begin{eqnarray}\label{near0}&
\int\limits_{-\delta}^{\delta}\Phi_h(x-y)v_0(y)dy=\int\limits_0^\delta(\Phi_h(x-y)-\Phi_h(x+y))v_0(y)dy\geq\cr
& \int\limits_0^\delta(\Phi_h(x-y)-\Phi_h(x+y))ky\,dy=\cr &
\int_{-\delta}^\delta \Phi_h(x-y)\kappa y\,dy = \kappa
\int_{-\delta-x}^{\delta-x} \Phi_h(z) (x+z)\,dz.
\end{eqnarray}
Combining \eqref{plato2} and \eqref{near0}, we obtain
\begin{eqnarray}\label{aboutv}&
v(x,h) \geq \kappa x
\int_{-\delta-x}^{\delta-x}\Phi_h(z)\,dz+H\int_{\delta-x}^{\delta+x}\Phi_h(z)\,dz
+ \kappa \int_{-\delta-x}^{\delta-x}z\Phi_h(z)\,dz -\cr &
C(\alpha)hxL^{-1-2\alpha}\|v_0\|_{L^\infty}.
\end{eqnarray}

Now we split the proof into several parts according to the regions
being considered. \\
I. Estimate for $0 \leq x \leq \delta=H/\kappa.$ Observe that for
$0\leq x\leq\delta/2,$ the contribution of the second and third
integrals in \eqref{aboutv} is positive. Indeed, it is equal to \[
\int_{\delta-x}^{\delta +x}(H - \kappa z)\Phi_h(z)\,dz, \] which
is positive due to monotonicity of $\Phi_h(z)$ and equality $H =
\kappa \delta.$ Thus, in this interval of $x$ we simply estimate
$v$ by dropping the combined contribution of the second and the
third integrals:
\begin{eqnarray}\nonumber &
v(x,h)\geq\kappa x
\int_{-\delta}^{\delta/2}\Phi_h(z)\,dz-C(\alpha)hxL^{-1-2\alpha}\|v_0\|_{L^\infty}\geq\cr
& \kappa
x(1-C(\alpha)h\delta^{-2\alpha})-C(\alpha)xhL^{-1-2\alpha}\|v_0\|_{L^\infty}
\geq \kappa x(1-C(\alpha)h\delta^{-2\alpha}).
\end{eqnarray}
Here we also decreased the interval of integration and used
\eqref{Phi}.

For $\delta/2\leq x\leq\delta$ we combine together the first and
the second integral and notice that $H=\kappa\delta\geq\kappa x$:
\begin{eqnarray}\nonumber &
v(x,h)\geq\kappa x \int_{-\delta-x}^{\delta+x}\Phi_h(z)\,dz+\kappa
\int_{-\delta-x}^{\delta-x}z\Phi_h(z)\,dz-C(\alpha)hxL^{-1-2\alpha}\|v_0\|_{L^\infty}\geq\cr
& \kappa x \int_{-\delta}^{\delta}\Phi_h(z)\,dz-\kappa
\int_{0}^{2\delta}z\Phi_h(z)\,dz-C(\alpha)hxL^{-1-2\alpha}\|v_0\|_{L^\infty}\geq\cr
& \kappa x(1-C(\alpha)h\delta^{-2\alpha})-C(\alpha)\kappa
h\delta^{1-2\alpha}-C(\alpha)xhL^{-1-2\alpha}\|v_0\|_{L^\infty}\geq
\kappa x(1-C(\alpha)h\delta^{-2\alpha}).
\end{eqnarray}
Here we again used \eqref{Phi} and \eqref{condstep}. Combining the
estimates together we have
\begin{equation}\label{smallxfin}
v(x,h)\geq\kappa x(1-C(\alpha)h\delta^{-2\alpha})
\end{equation}
for $0\leq x\leq\delta$.  \\

II. Estimate for $L-a \leq x \leq L$ case. The estimate is
virtually identical to the first case due to symmetry; $\delta$
has to be replaced by $a$. Thus, (we recall that
$\delta=H\kappa^{-1}\leq a$ by assumption of the lemma)
\begin{equation}\label{largexfin}
v(x,h) \geq \frac{H}{a} (L-x) (1-C(\alpha)h
a^{-2\alpha})\geq\frac{H}{a} (L-x) (1-C(\alpha)h
\delta^{-2\alpha}),
\end{equation}
for $L-a \leq x \leq L$.  \\

III. Estimate for $\delta \leq x \leq L/2.$ Here estimates are
similar to the first case. In the last term in \eqref{aboutv} we
will just estimate $x$ by $L$. Furthermore, observe that it
follows from \eqref{near0} and monotonicity property of $\Phi$
that the sum of the first and the third integrals in
\eqref{aboutv} is
positive. 
For $2\delta\leq x\leq L/2$ we ignore the positive combined
contribution of the first and the third integrals:
\begin{eqnarray}\nonumber &
v(x,h)\geq H
\int_{-\delta}^{\delta}\Phi_h(z)\,dz-C(\alpha)hL^{-2\alpha}\|v_0\|_{L^\infty}\geq\cr
&
H(1-C(\alpha)h\delta^{-2\alpha})-C(\alpha)hL^{-2\alpha}\|v_0\|_{L^\infty}
\geq H(1-C(\alpha)h\delta^{-2\alpha}).
\end{eqnarray}

For $\delta\leq x\leq2\delta$ we combine the first and the second
integrals and take into account that $\kappa x\geq\kappa\delta=H$:
\begin{eqnarray}\nonumber &
v(x,h)\geq H \int_{-\delta-x}^{\delta+x}\Phi_h(z)\,dz+\kappa
\int_{-\delta-x}^{\delta-x}z\Phi_h(z)\,dz-C(\alpha)hL^{-2\alpha}\|v_0\|_{L^\infty}\geq\cr
& H\int_{-\delta}^{\delta}\Phi_h(z)\,dz-\kappa
\int_{0}^{3\delta}z\Phi_h(z)\,dz-C(\alpha)hL^{-2\alpha}\|v_0\|_{L^\infty}\geq\cr
& H(1-C(\alpha)h\delta^{-2\alpha})-C(\alpha)\kappa
h\delta^{1-2\alpha}-C(\alpha)hL^{-2\alpha}\|v_0\|_{L^\infty}\geq
H(1-C(\alpha)h\delta^{-2\alpha}).
\end{eqnarray}
Combining the estimates together we obtain
\begin{equation}\label{platoest2}
v(x,h) \geq H (1 - C(\alpha)h \delta^{-2\alpha}).
\end{equation}

IV. Estimate for $L/2 \leq x \leq L-a.$ By symmetry we obtain
\begin{eqnarray}\label{platoest3}&
v(x,h) \geq H (1 - C(\alpha)ha^{-2\alpha})\geq H (1 -
C(\alpha)h\delta^{-2\alpha}).
\end{eqnarray}
Together, \eqref{platoest3}, \eqref{platoest2}, \eqref{largexfin} and
\eqref{smallxfin} complete the proof.
\end{proof}

Combining Proposition~\ref{timesplit} and Lemmas~\ref{burstep} and
\ref{disstep}, we obtain
\begin{theorem}\label{keystep}
Assume that the initial data $u_0(x)$ is $2L$-periodic, odd,
$u_0(L)=0,$ and $u_0(x) \geq \varphi(\kappa,H,a,x).$ Suppose that
\eqref{condstep} holds with $v_0$ replaced by $u_0.$ Assume also
that the solution $u(x,t)$ of the equation \eqref{bur1} with
initial data $u_0(x)$ satisfies $\|u(x,t)\|_{C^3} \leq N$ for $0
\leq t \leq T.$ Then for every $h \leq h_0(\alpha,N)$ small
enough, we have for $0 \leq x \leq L$
\begin{equation}\label{forest}
u(x,h) \geq \varphi(\tilde{\kappa},
\tilde{H},a+h\|u_0\|_{L^\infty},x),
\end{equation}
where \begin{equation}\label{newkap}
\tilde{\kappa}=\kappa(1-C(\alpha)\kappa^{2\alpha}
H^{-2\alpha}h)(1-\kappa h)^{-1} - C(\alpha, N)h^2
\end{equation}
and \begin{equation}\label{newH} \tilde{H}=
H(1-C(\alpha)\kappa^{2\alpha} H^{-2\alpha}h)-C(\alpha,N)h^2.
\end{equation}
\end{theorem}
\begin{proof}
We can clearly assume that $\kappa h  \leq 1/2;$ in view of our
assumptions on $u_0,$ $h \leq 1/2N$ is sufficient for that. Then
Lemmas~\ref{burstep} and \ref{disstep} together ensure that the
time splitting solution $v(x,h)$ of \eqref{weq} and \eqref{veq}
satisfies for $0\leq x \leq L$
\begin{equation}\label{vbound}
v(x,h) \geq \varphi(\kappa(1-C(\alpha)\kappa^{2\alpha}
H^{-2\alpha}h)(1-\kappa h)^{-1}, H(1-C(\alpha)\kappa^{2\alpha}
H^{-2\alpha}h), a+\|u_0\|_{L^\infty}h,x).
\end{equation}
Furthermore, Proposition~\ref{timesplit} allows us to pass from
the lower bound on $v(x,h)$ to lower bound on $u(x,h),$ leading to
\eqref{forest}, \eqref{newkap}, \eqref{newH}.
\end{proof}

From Theorem~\ref{keystep}, we immediately infer
\begin{corollary}\label{diffsys}
Under assumptions of the previous theorem and the additional
assumption stated below, for all $h$ small enough we have for $0
\leq x \leq L$ and $0 \leq nh \leq T$
\begin{equation}\label{unest}
u(x,nh) \geq \varphi(\kappa_n, H_n, a_n ,x).
\end{equation}
Here
\begin{equation}\label{kaprec}
\kappa_n = \kappa_{n-1} ( 1-C(\alpha)\kappa_{n-1}^{2\alpha}
H^{-2\alpha}_{n-1}h)(1-\kappa_{n-1} h)^{-1} - C(\alpha, N)h^2,
\end{equation}
\begin{equation}\label{Hrec2}
H_n = H_{n-1}(1-C(\alpha)\kappa_{n-1}^{2\alpha} H^{-2\alpha}_{n-1}h)
-C(\alpha,N)h^2,
\end{equation}
and
\begin{equation}\label{arec}
a_n=a+nh\|u_0\|_{L^\infty}.
\end{equation}
The corollary only holds assuming that for every $n,$ we have
\begin{equation}\label{Lcon}
H_n\kappa_n^{-1}\leq a_n,\ \ L\geq 4a_n,\ \
L^{-2\alpha}\|u_0\|_{L^\infty} \leq 4H_n a_n^{-2\alpha}.
\end{equation}
\end{corollary}

To study \eqref{kaprec} and \eqref{Hrec2}, we introduce the
following system of differential equations:
\begin{equation}\label{diffeqsys}
\kappa' = \kappa^2 - C(\alpha) \kappa^{1+2\alpha}
H^{-2\alpha};\,\,\,H'=-C(\alpha)\kappa^{2\alpha} H^{1-2\alpha}.
\end{equation}
\begin{lemma}\label{deapp}
Assume that $[0,T]$ is an interval on which the solutions of the
system \eqref{diffeqsys} satisfy $|\kappa (t)| \leq 2N,$
$0<H_1(\alpha) \leq H(t) \leq H_0(\alpha).$ Then for every
$\epsilon>0,$ there exists $h_0(\alpha,N,\epsilon)>0$ such that if
$h<h_0,$ then $\kappa_n$ and $H_n$ defined by \eqref{kaprec} and
\eqref{Hrec2} satisfy $|\kappa_n - \kappa(nh)|<\epsilon,$ $|H_n
-H(nh)| < \epsilon$ for every $n \leq [T/h].$
\end{lemma}
\begin{proof}
This is a standard result on approximation of differential
equations by a finite difference scheme. Observe that the
assumptions on $\kappa(t)$ and $H(t)$ also imply upper bounds on
$\kappa'(t),$ $\kappa''(t),$ $H'(t)$ and $H''(t)$ by a certain
constant depending only on $N$ and $\alpha.$ The result can be
proved comparing the solutions step-by-step inductively. Each step
produces an error not exceeding $C_1(\alpha,N)h^2,$ and the total
error over $[T/h]$ steps is estimated by $C_1(\alpha,N)h.$
Choosing $h_0(\alpha,N,\epsilon)$ sufficiently small completes the
proof.
\end{proof}
The final ingredient we need is the following lemma on the
behavior of solutions of the system \eqref{diffeqsys}.
\begin{lemma}\label{conlaw}
Assume that the initial data for the system \eqref{diffeqsys}
satisfy
\begin{equation}\label{iccon} H_0^{2\alpha} \kappa_0^{1-2\alpha}
\geq C(\alpha)/(1-2\alpha). \end{equation} Then on every interval
$[0,T]$ on which the solution makes sense (that is, $\kappa(t)$
bounded), the function $H(t)^{2\alpha} \kappa^{1-2\alpha}(t)$ is
non-decreasing.
\end{lemma}
\begin{proof}
A direct computation shows that
\[ \left( H(t)^{2\alpha} \kappa^{1-2\alpha}(t) \right)' = (1-2\alpha)
\kappa(t) \left(H(t)^{2\alpha} \kappa(t)^{1-2\alpha} -
\frac{C(\alpha)}{1-2\alpha}\right). \]
\end{proof}
Now we are ready to complete the blow up construction.
\begin{proof}[Proof of Theorem~\ref{thm3}]
Set $\kappa_0$ to be large enough, in particular
\begin{equation}\label{Hzero}
\kappa_0 = \left(\frac{3C(\alpha)}{1-2\alpha}
\right)^{\frac{1}{1-2\alpha}}
\end{equation}
will do. Set $H_0=1$, $a=\kappa_0^{-1}$,
$T(\alpha)=\frac{3}{2\kappa_0}$. Choose $L$ so that
\begin{equation}\label{Lchoice}
L \geq 16a.
\end{equation}
The initial data $u_0(x)$ will be a smooth, odd, $2L-$periodic
function satisfying $u_0(L)=0$ and $u_0(x) \geq
\varphi(\kappa_0,H_0,a,x).$ We will also assume
$\|u_0\|_{L^\infty} \leq 2H_0.$ Observe that $H_0$ and $\kappa_0$
are chosen so that in particular the condition \eqref{iccon} is
satisfied.
From
\eqref{diffeqsys} and Lemma~\ref{conlaw} it follows that
\begin{equation}\label{kappabu}
\kappa' = \kappa^2 -C(\alpha)\kappa^{1+2\alpha}H^{-2\alpha} \geq
\frac23\kappa^2.
\end{equation}
This implies $\kappa (t) \geq \frac{1}{\kappa_0^{-1}-\frac23 t}.$ In
particular, there exists $t_0 < T(\alpha)$ such that
$\kappa(t_0)=2N$ for the first time. Note that due to
\eqref{kappabu}, for $0 \leq t \leq t_0$ we have
\begin{equation}\label{kup}
 \kappa(t) \leq \frac{1}{\frac23(t_0-t)+\frac{1}{2N}} \leq
\frac{3/2}{(t_0-t)}.
\end{equation}
Rewrite the equation for $H(t)$ as
\begin{equation}\label{Halpha}
(H^{2\alpha})' = -2C(\alpha)\alpha \kappa^{2\alpha}.
\end{equation}
Using the estimate \eqref{kup} in \eqref{Halpha}, we get that for
any $0 \leq t \leq t_0,$
\[ H^{2\alpha}(t) \geq H_0^{2\alpha} - 2C(\alpha)\alpha
\int_0^{t_0}\kappa^{2\alpha}(s)\,ds \geq H_0^{2\alpha}(1-\alpha).
\] We used the fact that $H_0=1,$ $t_0 < T(\alpha) = \frac{3}{2\kappa_0}$
and \eqref{Hzero}. Now we can apply Lemma~\ref{deapp} on the
interval $[0,t_0].$ Choosing $\epsilon$ and $h$ sufficiently
small, we find that for $0 \leq nh \leq t_0,$ $\kappa_n \geq 1$
and $H_n \geq (1-\alpha)^{1/2\alpha}H_0\geq H_0/2.$ Also,
evidently, $a_n \leq a+2H_0 T(\alpha)=4a.$ This allows us to check
that the conditions \eqref{Lcon} hold on each step due to the
choice of $L$ \eqref{Lchoice}, justifying control of the true PDE
dynamics by the system \eqref{diffeqsys}.

From Lemma~\ref{deapp} and $\kappa(t_0)=2N$, we also see that,
given that $h$ is sufficiently small, $\kappa_{n_0} \geq 3N/2$ for
some $n_0$ such that $n_0h \leq t_0 < T(\alpha).$
Thus Corollary~\ref{diffsys} provides us with a lower bound
$u(0,n_0 h)=0,$ $u(x, n_0 h) \geq 3Nx/2$ for small enough $x.$
This contradicts our assumption that $\|u(x,t)\|_{C^3} \leq N$ for
$0 \leq t \leq T(\alpha),$ thus completing the proof.
\end{proof}

We obtained blow up in the case where period $2L$ was sufficiently
large (depending only on $\alpha$). However, examples of blow up
with arbitrary periodic data follow immediately from a scaling
argument. Indeed, assume $u(x,t)$ is a $2L-$periodic solution of
\eqref{bur1}. Then $u_1(x,t) = L^{-1+2\alpha} u(Lx,L^{2\alpha} t)$
is a $2-$periodic solution of the same equation. Thus a scaling
procedure allows to build blow up examples for any period.

\it Remark. \rm Formally we proved the blow up only in $C^3$
class. But since global regularity in $H^s$ class for
$s>\frac32-2\alpha$ provides global regularity in $C^\infty$ (see
Theorem~\ref{thm1}), we can conclude that we constructed a blow up
in $H^s$ class for every
$s>\frac32-2\alpha$. \\

\section{Global existence and regularity for rough initial data for the case
$\alpha=1/2$}\label{rough}

In this section we present some results on existence of regular
solution for $\alpha=1/2$ and rough initial data. More precisely,
we prove that the solution becomes smooth starting from any
initial data of the class $L^p,\ p>1$. It is natural that the
result can be obtained for the case $\alpha>1/2$ by more
traditional means. 
In the
present section we consider the case $\alpha=1/2$ only.

Consider the equation
\begin{equation}\label{eq}
u_t=uu_x-(-\Delta)^{1/2}u,\ \ u(x,0)=u_0(x),
\end{equation}
with $u_0\in L^p$ for some $p>1$. Let us look first at the
approximating equation
\begin{equation}\label{eq1}
u_t^N=u^Nu_x^N-(-\Delta)^{1/2}u^N,\ \ u^N(x,0)=u_0^N(x),
\end{equation}
where $u^N_0\in C^{\infty}$ and $\|u_0^N-u_0\|_{L^p}\to0$ as $N\to\infty$.
%
We need the following fact.
\begin{lemma}\label{lpnorms}
Assume that a smooth function $w(x,t)$ satisfies the equation
\eqref{eq} with smooth initial data $w_0(x).$ Then for every $1 <
p \leq \infty$ and every $t,$ we have $\|w(x,t)\|_{L^p} \leq
\|w_0(x)\|_{L^p}.$
\end{lemma}
This Lemma can be proven in the same way as a corresponding result
for the quasi-geostrophic equation, using the positivity of $\int
|w|^{p-2} w (-\Delta)^{\alpha}w\,dx.$ See \cite{Resnick} or
\cite{CC1} for more details.

We divide our proof of regularity into three steps.

{\it Step I.} Here we prove uniform (in $N$) estimates for the
$L^\infty$ norm.
Put
$$
M_N(t):=\|u^N(\cdot,t)\|_{L^\infty}.
$$
Fix $t\geq0$. Consider any point $x_0$ where $|u^N(x_0,t)|=M_N$.
Without loss of generality, we may assume that $x_0=0$ and
$u^N(0,t)=M_N.$ Then
\begin{equation}\label{equat}
u^N_t(0,t)=(-(-\Delta)^{1/2}u^N)(0,t)= \frac{1}{\pi}\int\limits_{-\infty}^{\infty}\frac{u^N(y,t)-M_N}{y^{2}}dy.
\end{equation}
Denote Lebesgue measure of a measurable set $S$ by $m(S).$ Since
by Lemma~\ref{lpnorms} we have
$$
\|u^N\|_{L^p}^p\leq C,
$$
we obtain that \[ m\left.\left(x \right| |u^N(x,t)| \geq M_N/2
\right) \leq C2^p M_N^{-p}.\] Then the right hand side of
\eqref{equat} does not exceed
$$
-M_N\int\limits_{L\geq|y|\geq C2^{p-1}/M_N^p}y^{-2}dy.
$$
Here $2L$ is the period. Then
\begin{equation}\label{infin1}
u^N_t(0,t)< -C_1 M_N^{p+1}+C_2 M_N.
\end{equation}
The same bound holds for any point $x_0$ where $M_N$ is attained
and by continuity in some neighborhoods of such points. So, we
have \eqref{infin1} in some open set $U_N$. Due to smoothness of
the approximating solution, away from $U_N$ we have
$$
\max\limits_{x\not\in U_N}|u^N(x,\tau)|< M_N(\tau)
$$
for every $\tau$ during some period of time $[t,t+\tau_N],$
$\tau_N>0.$ Thus we obtain that
\begin{equation}\label{infin}
\frac{d}{dt}M_N< -C_1 M_N^{p+1}+C_2 M_N.
\end{equation}
Solving equation \eqref{infin}, we get the uniform estimate
$$
M_N^{p}(t)\leq
\frac{e^{C_2pt}}{M_N^{-p}(0)+\frac{C_1}{C_2}(e^{C_2pt}-1)}\leq
\frac{C_2}{C_1(1-e^{-C_2pt})}.
$$
In particular,
\begin{equation}\label{linf}
t^{1/p}\|u^N\|_{L^\infty}\leq C,\ \ \ t\leq1.
\end{equation}

{\it Step II.} Here we obtain uniform in $N$ estimates on the
approximations $u^N$ that will imply smoothness of the solution.
We will use the construction similar to the one appearing in the
proof of Theorem~\ref{critcase}.

Clearly, it is sufficient to work with $t \leq 1.$ Let us define
\[ G(t) = {\rm inf}_{0 \leq \omega(x) \leq Ct^{-1/p}} \frac{\omega(x)}{x},
\]
where $C$ is as in \eqref{linf}. Observe that, since $\omega$ is
concave and increasing, the function $G(t)$ is equal to
$Ct^{-1/p}/\omega^{-1}(Ct^{-1/p}).$ Define also
\begin{equation}\label{deff}
F(t) =\left( \int\limits_0^t G(s)\,ds \right)^{-1}.
\end{equation}
 We claim that
solution $u^N(x,t)$ has modulus of continuity $\omega_{F(t)}$ for
every $t>0$ and every $N.$ Here $\omega$ is defined by
\eqref{modcon}. Let us fix an arbitrary $N>0.$ Since $u_0^N$ and
$u^N$ are both smooth and $F(t) \rightarrow \infty$ as $t
\rightarrow \infty,$ it follows that $u^N(x,t)$ has
$\omega_{F(t)}$ for all $t < t_0(N),$ $t_0(N)>0.$ By the argument
completely parallel to that of Lemma~\ref{scen11}, we can show
that if the modulus of continuity $\omega_{F(t)}$ is ever
violated, then there must exist $t_1 >0$ and $x \ne y$ such that
\[ u^N(x,t_1) - u^N(y,t_1) = \omega(F(t_1)|x-y|) \]
and $u^N(x,t)$ has $\omega_{F(t)}$ for any $t\leq t_1.$ Let us
denote $|x-y|$ by $\xi.$ Now consider \begin{eqnarray}
\label{rr11} \frac{\partial}{\partial t} \left. \left[
\frac{u^N(x,t)-u^N(y,t)}{\omega(F(t)\xi)} \right] \right|_{t=t_1}
= \frac{\partial_t (u^N(x,t)-u^N(y,t))|_{t=t_1}
\omega(F(t_1)\xi)}{\omega(F(t_1)|x-y|)^2} \\ \nonumber -
\frac{\omega(F(t_1)\xi)F'(t_1)\xi
\omega'(F(t_1)\xi)}{\omega(F(t_1)|x-y|)^2}.
\end{eqnarray}
It follows from the proof of Theorem~\ref{critcase} (see
\eqref{for1}, \eqref{for2}) that
\[ \left. \frac{d}{dt} \left( u^N(x,t)-u^N(y,t)
\right)\right|_{t=t_1} < -
\omega(F(t_1)\xi)\frac{d}{d\xi}\omega(F(t_1) \xi). \] Thus the
numerator on the right hand side of \eqref{rr11} is smaller than
\[ -\omega(F(t_1)\xi)^2 \omega'(F(t_1)\xi)F(t_1) -
\omega(F(t_1)\xi) F'(t_1)\omega'(F(t_1)\xi)\xi. \] The numerator
is strictly negative as far as
\begin{equation}\label{rr12}
-\frac{F'(t_1)}{F(t_1)^2} \leq
\frac{\omega(F(t_1)\xi)}{F(t_1)\xi}.
\end{equation}
Notice that by \eqref{linf}, we have
\[ \omega(F(t_1)\xi) =u^N(x,t_1)-u^N(y,t_1) \leq
2Ct_1^{-1/p}. \] Using the definition of the function $G(t),$ we
obtain the estimate
\[ \frac{\omega(F(t_1)\xi)}{F(t_1)\xi} \geq G(t_1). \]
Thus \eqref{rr12} is satisfied if
\[ \left(\frac1F \right)' \leq G(t), \]
which is correct by definition of $F.$ Therefore we obtain
\[ \frac{\partial}{\partial t} \left. \left[
\frac{u^N(x,t)-u^N(y,t)}{\omega(F(t)\xi)} \right] \right|_{t=t_1}
<0. \] Since $N$ was arbitrary, it follows that $u(x,t)$ has the
modulus of continuity $\omega_{F(t)}$ for any $t>0,$ and thus
\begin{equation}
\label{last} F(t)\|u^N\|_{W^1_\infty}\leq C,\ \ \ t\leq1.
\end{equation}
To obtain higher order regularity of the solution we apply
arguments from the proof of Theorem~\ref{existcinfgal}. We start
with $s=1,\ q=1$. We can repeat the proof step by step. The only
difference is that now we will use the estimate \eqref{last}
instead of uniform bound for the norm $\|u^N\|_1$. Finally, we
obtain the estimates
\begin{equation}\label{last1}
F_n(t)\|u^N(\cdot,t)\|_{1+\frac{n}{2}}\leq C_n,\ \ \ n\geq1,\ \ \ t\leq1,
\end{equation}
with some functions $F_n$ which can be calculated inductively.
Now, we can choose a subsequence $N_j$ (cf. proof of
Theorem~\ref{existcinf}) such that $u^{N_j}\to u$ as
$N_j\to\infty$ and function $u$ satisfies differential equation
\eqref{eq} as well as the bounds \eqref{linf}, \eqref{last},
\eqref{last1} on $(0,1]$.

{\it Step III.} Here we prove that the function $u$ can be chosen
to satisfy the initial condition.

\begin{lemma}\label{converg}
Assume that $p\in(1,\infty)$. Then $\|u(\cdot,t)-u_0(\cdot)\|_{L^p}\to0$ as $t\to0$.
\end{lemma}
\begin{proof}
Let $\varphi(x)$ be an arbitrary $C^{\infty}$ function. Put
$$
g^N(t,\varphi):=(u^N,\varphi)=\int u^N(x,t)\varphi(x)dx.
$$
Obviously, $g^N(\cdot,\varphi)\in C([0,1])$. We will use the
estimate \eqref{estim}:
\begin{equation}\label{estimnew}
\int_0^1 |g_t^N|^{1+\delta}\,dt \leq
C\left(\int_0^1\|u^N\|_{L^2}^{2+2\delta}
\|\varphi\|^{1+\delta}_{W^1_\infty}\,dt + \int_0^1
\|u^N\|^{1+\delta}_{L^2} \|\varphi\|_{2\alpha}^{1+\delta}\,dt
\right),
\end{equation}
which holds for any $\delta>0$.
Put $\delta:=1$ if $p\geq2$ and $\delta:=(p-1)/(2-p)$ if $1<p<2$.
Due to \eqref{linf} we obtain
\begin{equation}\label{l2}
\|u^N\|_{L^2}^2\leq\|u^N\|_{L^\infty}^{2-p}\|u^N\|_{L^p}^p\leq
Ct^{-\frac{2-p}{p}},\ \ \ t\leq1.
\end{equation}
Substituting \eqref{l2} into \eqref{estimnew} we see that
$\|g^N_t(\cdot,\varphi)\|_{L^{1+\delta}}\leq C(\varphi)$. By the
same argument as used in the proof of Theorem~\ref{existcinf} we
conclude that there exists a subsequence $u^{N_j}$ such that
for any $\varphi\in L^{p'}$ the sequence $(u^{N_j},\varphi)$ tends
to $(u,\varphi)$ uniformly on $[0,1]$.

Next,
\begin{equation}\label{fin}
|(u-u_0,\varphi)|\leq
|(u-u^{N_j},\varphi)|+|(u^{N_j}-u^{N_j}_0,\varphi)|+|(u^{N_j}_0-u_0,\varphi)|.
\end{equation}
The first and the third terms in the right hand side of the
\eqref{fin} can be made small uniformly on $[0,1]$ by choosing
sufficiently large $N_j$. The second term tends to zero as $t\to0$
for every fixed $N_j$. Thus $u(\cdot,t)$ converges to $u_0(\cdot)$
as $t\to0$ weakly in $L^p$. In particular,
\begin{equation}\label{upper}
\|u_0(\cdot)\|_{L^p}\leq \liminf_{t\to0}\|u(\cdot,t)\|_{L^p}.
\end{equation}

Due to monotonicity property of Lemma~\ref{lpnorms} we have
\begin{equation}\label{below}
\|u_0(\cdot)\|_{L^p}\geq \limsup_{t\to0}\|u(\cdot,t)\|_{L^p}.
\end{equation}
Thus $\|u_0(\cdot)\|_{L^p}=\lim_{t\to0}\|u(\cdot,t)\|_{L^p}$. Now,
it follows from uniform concavity of the space $L^p$,
$p\in(1,\infty)$, that weak convergence and convergence of norms
imply convergence in the norm sense.
\end{proof}

Let us combine the results in the following theorem.
\begin{theorem}\label{thmrough}
Let $u_0\in L^p$ for some $p\in(1,\infty)$. Then there exists a solution $u(x,t)$
of the equation \eqref{eq} such that $u$ is real analytic for $t>0$,
\begin{equation}
\|u(\cdot,t)-u_0(\cdot)\|_{L^p}\to0\ \ \ \hbox{as}\ \ t\to0;
\end{equation}
\begin{equation}
t^{1/p}\|u(\cdot,t)\|_{L^\infty}\leq C(\|u_0\|_{L^p}),\ \ \
0<t\leq1;
\end{equation}
\begin{equation}
F(t)^{-1}\|u(\cdot,t)\|_{W^1_\infty}\leq C(\|u_0\|_{L^p}),\ \ \
0<t\leq1;
\end{equation}
Here $F$ is defined in \eqref{deff}.
\end{theorem}

\noindent {\it Remark} 1. If $u_0\in H^s$ for some $s>1/2$ then
$u$ converges to $u_0$ in the $H^s$ norm as well. However, the
question whether we have convergence in $H^r$ norm if $u_0\in
H^r$, $0<r<1/2$, is still open. For the case $r=1/2$ the answer is
positive (see Section 6).

\noindent \it Remark \rm 2. Another interesting open question is
the uniqueness of the solution from Theorem~\ref{thmrough}. Due to
the highly singular nature of estimates as $t$ approaches zero,
the usual uniqueness argument based on some sort of Gronwall
inequality does not seem to go through.

\section{The critical Sobolev space}\label{critical}

Here we show that the results of Theorems~\ref{subcritical},
\ref{thm1}, \ref{thm2}, \ref{thm3}, \ref{existcinf}, \ref{critcase} and
\ref{analcrit} hold for $s=3/2-2\alpha$, as well.

Assume that $u_0\in H^s,\ s\geq q=3/2-2\alpha$, $1>\alpha>0$. We
introduce the following Hilbert spaces of periodic functions. Let
$\varphi: [0,\infty)\to[1,\infty)$ be an unbounded increasing
function. Then $H^{s,\varphi}$ consists of periodic functions
$f\in L^2$ such that its Fourier coefficients satisfy
\begin{equation}
\|f\|_{H^{s,\varphi}}^2:=\sum\limits_n|n|^{2s}(\varphi(|n|))^2|\hat{f}(n)|^2<\infty.
\end{equation}
Note that $u_0\in H^{s,\varphi}$ for some function $\varphi$.
Without loss of generality we may assume, in addition, that
$\varphi\in C^{\infty}$ and
\begin{equation}\label{phiest}
\varphi'(x)\leq Cx^{-1}\varphi(x)
\end{equation}
for some constant $C$. It follows from \eqref{phiest} that
\begin{equation}\label{phiest1}
\varphi(2x)\leq 2^C\varphi(x).
\end{equation}
We start from Galerkin approximations. Consider the sum arising
from the nonlinear term when estimating the $H^s$ norm of the
solution:
\[S:=\sum\limits_{a+b+k=0,|a|,|b|,|k|\leq N}k|k|^{2s}(\varphi(|k|))^2\hat{u}^N(a)\hat{u}^N(b)\hat{u}^N(k).\]
In what follows, for the sake of brevity, we will omit mentioning
restrictions $|a|,|b|,|k|\leq N$ in notation for the sums; all
sums will be taken with this restriction. Observe that (cf.
\eqref{sums11})
\begin{equation}\label{sumsnew}
|S|\leq 6\sum\limits_{k+a+b=0,|a|\leq|b|\leq|k|}
|k|k|^{2s}(\varphi(|k|))^2+a|a|^{2s}(\varphi(|a|))^2+b|b|^{2s}(\varphi(|b|))^2||\hat{u}^N(k)\hat{u}^N(a)\hat{u}^N(b)|.
\end{equation}
Recall that under conditions $|a|\leq|b|\leq|k|,\ \ a+b+k=0$, we
have $|a|\leq|k|/2,\ |b|\geq|k|/2$. Next, due to \eqref{phiest}
and \eqref{phiest1} we estimate
\begin{equation}\label{teches11}
\begin{split}
&
|k|k|^{2s}(\varphi(|k|))^2+a|a|^{2s}(\varphi(|a|))^2+b|b|^{2s}(\varphi(|b|))^2|=
\cr &
|b(|b|^{2s}(\varphi(|b|))^2-|b+a|^{2s}(\varphi(|b+a|))^2)+a(|a|^{2s}(\varphi(|a|))^2-|k|^{2s}(\varphi(|k|))^2|\leq
\cr & C|a||k|^{2s}(\varphi(|k|))^2\leq
C|a||b|^{s}\varphi(|b|)|k|^{s}\varphi(|k|).
\end{split}
\end{equation}
Fix $M>0$ to be specified later. Notice that sum over $|k| \leq M$
in \eqref{sumsnew} can be bounded by a constant $C(M).$ Splitting
summation in $a$ over dyadic shells scaled with $|k|,$ define
\[ S_1(l) = \sum\limits_{k+a+b=0,|b|\leq|k|,|k|\geq M,|a|\in[2^{-l-1}|k|,2^{-l}|k|]}
|a|^{1-2\alpha}|b|^{s+\alpha}\varphi(|b|)|k|^{s+\alpha}\varphi(|k|)|\hat{u}^N(k)\hat{u}^N(a)\hat{u}^N(b)|.
\]
Then due to \eqref{teches11} and the relationship between $a,$ $b$
and $k$ in the summation for $S$ we have
\begin{equation}\label{abouts}
|S|\leq C\sum\limits_{l=1}^\infty 2^{-2l\alpha}S_1(l)+C(M).
\end{equation}
Think of $S_1(l)$ as a quadratic form in $\hat{u}^N(k)$ and
$\hat{u}^N(b).$ Then applying Schur test to each $S_1(l)$ we
obtain
\begin{equation}
\begin{split}
&
S_1(l)\leq\|u^N\|_{H^{s+\alpha,\varphi}}^2\cdot\sup\limits_{|k|\geq
M}\sum\limits_{|a|\in[2^{-l-1}|k|,2^{-l}|k|]}|a|^{1-2\alpha}|\hat{u}^N(a)|\leq
\cr &
C\|u^N\|_{H^{s+\alpha,\varphi}}^2\|u^N\|_{H^{q,\varphi}}(\varphi(2^{-l}M))^{-1}.
\end{split}
\end{equation}
Next, note that
\begin{equation}\label{abouts1}
\begin{split}
& \sum\limits_{l=1}^\infty
2^{-2l\alpha}S_1(l)=\sum\limits_{l=1}^{l_0}
2^{-2l\alpha}S_1(l)+\sum\limits_{l=l_0}^\infty
2^{-2l\alpha}S_1(l)\leq \cr &
C\|u^N\|_{H^{s+\alpha,\varphi}}^2\|u^N\|_{H^{q,\varphi}}
\left(\frac{1}{1-2^{-2\alpha}}(\varphi(2^{-l_0}M))^{-1}+\frac{2^{-2l_0\alpha}}{1-2^{-2\alpha}}\right).
\end{split}
\end{equation}
Given $\epsilon>0$, we can choose, first, sufficiently large $l_0$
and then sufficiently large $M$ to obtain from \eqref{abouts},
\eqref{abouts1} and unboundedness of $\varphi$
\begin{equation}\label{aboutsnew}
|S|\leq
C\epsilon\|u^N\|_{H^{s+\alpha,\varphi}}^2\|u^N\|_{H^{q,\varphi}}+C(M(\epsilon)).
\end{equation}
It follows from \eqref{gal} and \eqref{aboutsnew} that
\begin{equation}
\frac{d}{dt}\|u^N\|_{H^{s,\varphi}}^2\leq
(C\epsilon\|u^N\|_{H^{q,\varphi}}-1)\|u^N\|_{H^{s+\alpha,\varphi}}^2+C(\epsilon),\
\ s\geq q=3/2-2\alpha,\ \alpha>0.
\end{equation}
Using this estimate and the same arguments as before we can extend
the results of Theorems~\ref{subcritical}, \ref{thm1}, \ref{thm2}, \ref{thm3},
\ref{existcinf}, \ref{critcase} and \ref{analcrit} to the case
$s=3/2-2\alpha,\ \alpha>0$. Here we formulate them for convenience
of future references.
\begin{theorem}\label{subcriticalnew}
Assume that $\alpha>1/2,$ and the initial data $u_0(x) \in H^s,$
$s\geq3/2-2\alpha,$ $s \geq 0.$ Then there exists a global
solution of the equation \eqref{bur1} $u(x,t)$ which belongs to
$C([0,\infty),H^s)$ and is real analytic in $x$ for $t>0$.
\end{theorem}
\begin{theorem}\label{critcasenew}
Assume $\alpha=1/2,$ and $u_0 \in H^s,$ $s\geq1/2.$ Then there
exists a global solution $u(x,t)$ of \eqref{bur1} which belongs to
$C([0,\infty),H^s)$ and is real analytic in $x$ for $t>0$. If $v$
is another weak solution of \eqref{bur1} with initial data $u_0$
such that $v\in C([0,T],L^2)\cap
L^{3/2\delta}([0,T],H^\delta)$ with some $\delta\in(1/2,1]$, then $v$
coincides with $u$ on $[0,T]$.
\end{theorem}
\begin{theorem}\label{thm1new}
Assume that $0<\alpha<1/2,$ and the initial data $u_0(x) \in H^s,$
$s\geq3/2-2\alpha$. Then there exists $T=T(\alpha,u_0)>0$
such that there exists a weak solution of the equation
\eqref{bur1} $u(x,t) \in C([0,T],H^s) \cap
L^2([0,T],H^{s+\alpha}).$ Moreover, $u(x,t) \in C^{\infty}$ for
any $0<t<T.$ If $v$ is another weak solution of \eqref{bur1} with
initial data $u_0$ such that $v\in C([0,T],L^2)\cap
L^{3/2\delta}([0,T],H^{\delta})$ with some $\delta\in(1/2,1]$, then $v$ coincides
with $u$.
\end{theorem}
\begin{theorem}\label{thm3new}
Assume that $0<\alpha<1/2.$ Then there exists smooth periodic
initial data $u_0(x)$ such that the solution $u(x,t)$ of
\eqref{bur1} blows up in $H^s$ for each $s\geq\frac32 -2\alpha$ in
a finite time.
\end{theorem}
\begin{theorem}\label{existcinfnew}
Assume that $s \geq 3/2-2\alpha$, $s\geq0$, $\alpha>0$, and $u_0
\in H^s.$ Then there exists $T=T(\alpha,u_0)>0$ and a solution $u(x,t)$ of \eqref{bur1}
such that
\begin{equation}\label{ws2new}
u \in L^2([0,T], H^{s+\alpha}) \cap C([0,T],H^s);
\end{equation}
\begin{equation}\label{cs2new}
t^{n/2}u  \in  C((0,T],H^{s+n \alpha})\cap L^\infty([0,T],H^{s+n
\alpha})
\end{equation}
for every $n>0.$
\end{theorem}
{\it Remark.} If $s>3/2-2\alpha$ then
$T(\alpha,u_0)=T(\alpha,\|u_0\|_s)$. If $s=3/2-2\alpha$ then
$u_0\in H^{s,\varphi}$ for some function $\varphi$ described at
the beginning of the section and
$T(\alpha,u_0)=T(\alpha,\varphi,\|u_0\|_{H^{s,\varphi}})$.

{\bf Acknowledgement.} AK and RS have been supported in part by
NSF grant DMS-0314129. FN has been supported in part by NSF grant
DMS-0501067. AK thanks Igor Popov for useful discussions.


\begin{thebibliography}{99}

\bibitem{ADV} N.~Alibaud, J.~Droniou and J.~Vovelle, \it Occurence and
non-appearance of shocks in fractal Burgers equations, \rm
preprint

\bibitem{BM} A.~Bertozzi and A.~Majda, \it Vorticity and
Incompressible Flow, \rm Cambridge University Press, 2002

\bibitem{Beteman} H.~Beteman, \it Some recent researches of the
motion of fluid, \rm Monthly Weather Rev. {\bf 43} (1915),
163--170

\bibitem{CafVas} L. Caffarelli, A. Vasseur,
\it Drift diffusion equations with fractional diffusion and the
quasi-geostrophic equation, \rm preprint

\bibitem{CF} P.~Constantin and C.~Foias, \it Navier-Stokes
Equations, \rm The University of Chicago Press, 1989

\bibitem{CMT} P.~Constantin, A.~Majda and E.~Tabak, \it Formation
of strong fronts in the 2D quasi-geostrophic thermal active
scalar, \rm Nonlinearity, {\bf 7} (1994), 1495--1533

\bibitem{CCW} P.~Constantin, D.~Cordoba and J.~Wu, \it On the critical dissipative
quasi-geostrophic equation, \rm Dedicated to Professors Ciprian
Foias and Roger Temam (Bloomington, IN, 2000). Indiana Univ. Math.
J. {\bf 50} (2001), 97--107

\bibitem{Cord1} D.~Cordoba, \it Nonexistence of simple hyperbolic
blow up for the quasi-geostrophic equation, \rm Ann. of Math. {\bf
148} (1998), 1135--1152

\bibitem{CC} A.~Cordoba and D.~Cordoba, \it A maximum principle
applied to quasi-geostrophic equations, \rm Commun. Math. Phys.
{\bf 249} (2004), 511--528

\bibitem{CC1} A.~Cordoba and D.~Cordoba, \it A pointwise estimate for
fractionary derivatives with applications to partial differential
equations, \rm Proc. Natl. Acad. Sci. USA {\bf 100} (2003),
15316--15317

\bibitem{Feller} W.~Feller, \it Introduction to Probability Theory and
Its Applications, \rm Vol. {\bf 2}, Wiley, 1971

\bibitem{Forsyth} A.R.~Forsyth, \it Theory of Differential
Equations, \rm Vol. {\bf 6}, Cambridge University Press, 1906

\bibitem{KP1} N.~Katz and N.~Pavlovi\'c, \it A cheap Caffarelli-Kohn-Nirenberg
inequality for the Navier-Stokes equation with hyper-dissipation,
\rm Geom. Funct. Anal. {\bf 12} (2002), 355--379

\bibitem{KNV} A.~Kiselev, F.~Nazarov and A.~Volberg, \it Global well-posedness for
the critical $2D$ dissipative quasi-geostrophic equation, \rm
Inventiones Math. {\bf 167} (2007) 445--453

\bibitem{Matalon} M.~Matalon, \it Intrinsic flame instabilities in
premixed and non-premixed combustion, \rm Annu Rev. Fluid Mech.
{\bf 39 }(2007), 163--191

\bibitem{Resnick} S.~Resnick, \it Dynamical problems in nonlinear
advective partial differential equations, \rm Ph.D. Thesis,
University of Chicago, 1995

\bibitem{T3} M.~Taylor, \it Partial Differential Equations III:
Nonlinear Equations, \rm Springer-Verlag, New York, 1997

\bibitem{Wu} J.~Wu, \it The quasi-geostrophic equation
and its two regularizations, \rm Comm. Partial Differential
Equations {\bf 27} (2002), 1161--1181

\end{thebibliography}
\end{document}